\documentclass[]{sn-jnl}

\usepackage{graphicx}%
\usepackage{multirow}%
\usepackage{amsmath,amssymb,amsfonts}%
\usepackage{amsthm}%
\usepackage{mathrsfs}%
\usepackage[title]{appendix}%
\usepackage{xcolor}%
\usepackage{textcomp}%
\usepackage{manyfoot}%
\usepackage{booktabs}%
\usepackage{algpseudocode}%
\usepackage{listings}%

\usepackage{graphics}
\usepackage{tikz}
\usepackage{pgfplots}
\usepackage{booktabs}
\usepackage{pgfplotstable}
\pgfplotsset{compat = 1.3}
\pgfplotstableset{
        every head row/.style={before row=\toprule, after row=\midrule},
        every last row/.style={after row=\bottomrule},
}
\usetikzlibrary{fillbetween}
\usepgfplotslibrary{fillbetween}
\usepackage{xparse}
\usepackage{bigints}

\usepackage{url}
\usepackage{hyperref}
\usepackage[nameinlink, capitalize, noabbrev]{cleveref}

\usepackage{mathtools}
\usepackage{subfigure}

\newcommand{\R}{\mathbb R}

\newcommand{\Prod}{p}
\newcommand{\dest}{d}

\newcommand{\diff}[1]{{\mathrm{d}{#1}}}

\newcommand{\ddt}[1]   {\frac{\partial{#1}}{\partial{t}}}
\newcommand{\ddx}[1]   {\frac{\partial{#1}}{\partial{x}}}
\newcommand{\ddy}[1]   {\frac{\partial{#1}}{\partial{y}}}

\newcommand{\ddo}[2]   { \frac{\diff #1}{\diff #2}}

\newcommand{\bu}{\mathbf{u}}

\newcommand{\bF}{\mathbf{F}}
\newcommand{\bG}{\mathbf{G}}
\newcommand{\bbS}{\mathbf{S}}

\newcommand{\bU}{\mathbf{U}}

\newcommand{\xip}{x_{i+1/2}}
\newcommand{\xin}{x_{i-1/2}}

\newcommand{\yjp}{y_{j+1/2}}
\newcommand{\yjn}{y_{j-1/2}}

\newcommand{\iip}{i+1/2}
\newcommand{\iin}{i-1/2}

\newcommand{\jjp}{j+1/2}
\newcommand{\jjn}{j-1/2}

\newcommand{\bS}{\mathbf{S}}

\usepackage{bbm}

\def\R{\mathbb{R}}

\def\dt{\Delta t}
\definecolor{darkspringgreen}{rgb}{0., 0.55, 0.3}
\definecolor{dartmouthgreen}{rgb}{0.05, 0.5, 0.06}
\definecolor{etonblue}{rgb}{0.59, 0.78, 0.64}
\definecolor{airforceblue}{rgb}{0., 0.4, 0.66}
\definecolor{arylideyellow}{rgb}{0.91, 0.84, 0.42}
\definecolor{emerald}{rgb}{0.31, 0.78, 0.47}
\definecolor{uclagold}{rgb}{1.0, 0.7, 0.0}
\definecolor{cadmiumorange}{rgb}{0.93, 0.53, 0.18}

\newcommand{\uvec}[2][3]{\boldsymbol{#2\mkern-#1mu}\mkern#1mu}

\theoremstyle{thmstyleone}%
\theoremstyle{thmstyletwo}%
\newtheorem{remark}{Remark}%

\theoremstyle{thmstylethree}%

\numberwithin{equation}{section}

\begin{document}

\title[Article Title]{A high-order, fully well-balanced, unconditionally positivity-preserving finite volume framework for flood simulations}

\author[1]{\fnm{Mirco} \sur{Ciallella}}\email{mirco.ciallella@ensam.eu}

\author[2]{\fnm{Lorenzo} \sur{Micalizzi}}\email{lmicali@ncsu.edu}

\author[3]{\fnm{Victor} \sur{Michel-Dansac}}\email{victor.michel-dansac@inria.fr}

\author[4]{\fnm{Philipp} \sur{\"Offner}}\email{mail@philippoeffner.de}

\author[5]{\fnm{Davide} \sur{Torlo}}\email{davide.torlo@sissa.it}


\affil[1]{\orgname{\'Ecole Nationale Sup\'erieure d'Arts et M\'etiers, I2M}, \city{Bordeaux}, \country{France}}

\affil[2]{\orgdiv{Department of Mathematics}, \orgname{North Carolina State University}, \city{Raleigh} \country{United States}}

\affil[3]{\orgname{Universit\'e de Strasbourg, CNRS, Inria, IRMA}, \city{Strasbourg}, \country{France}}

\affil[4]{\orgdiv{Institute of Mathematics}, \orgname{Johannes Gutenberg-University Mainz and TU Clausthal}, \city{Clausthal-Zellerfeld}, \country{Germany}}

\affil[5]{\orgdiv{SISSA mathLab}, \orgname{SISSA}, \city{Trieste}, \country{Italy}}


\abstract{In this work, we present a high-order finite volume framework for the numerical simulation of shallow water flows.
The method is designed to accurately capture complex dynamics inherent in shallow water systems, particularly suited for applications such as tsunami simulations.
The arbitrarily high-order framework ensures precise representation of flow behaviors,
crucial for simulating phenomena characterized by rapid changes and fine-scale features.
Thanks to an {\it ad-hoc} reformulation in terms of production-destruction terms,
the time integration ensures positivity preservation without any time-step restrictions,
a vital attribute for physical consistency, especially in scenarios where negative water depth
reconstructions could lead to unrealistic results.
In order to introduce the preservation of general steady equilibria dictated by the underlying balance law,
the high-order reconstruction and numerical flux are blended in a convex fashion with a well-balanced approximation,
which is able to provide exact preservation of both static and moving equilibria.
Through numerical experiments, we demonstrate the effectiveness and robustness of the proposed approach in capturing
the intricate dynamics of shallow water flows, while preserving key physical properties essential for flood simulations.}


\maketitle

\section{Introduction}

The Saint-Venant equations, also known as the shallow water (SW) equations, characterize the behavior of hydrostatic free surface waves influenced by gravity.
These nonlinear hyperbolic partial differential equations (PDEs) are applicable under the assumption of either very large wavelengths or very shallow depths.
They find extensive use across engineering domains, including river and estuarine hydrodynamics, urban flood management, and tsunami risk evaluation.
The numerical approximation of the SW equations remains a highly active area of research.
Numerous original methods have been developed across various contexts and settings:
finite volume~\cite{audusse2004fast,gallardo2007well,noelle2007high,diaz2013high,bollermann2013well,kurganov2002central,kurganov2018finite,cheng2016moving,xing2011advantage,berthon2016fully,michel2017well,michel2021two,ciallella2022arbitrary,ciallella2023arbitrary}, continuous and discontinuous finite element~\cite{https://doi.org/10.1002/fld.1650211009,https://doi.org/10.1002/fld.250,cite-keyi2,https://doi.org/10.1002/fld.1153,cite-key1,Scovazzi3,BEHZADI2020112662,xing2014exactly,mantri2021well,arpaia2022efficient,bender2024entropy,wintermeyer2015entropy}, residual distribution~\cite{ricchiuto2009stabilized,ricchiuto2007application,ricchiuto2015explicit,ricchiuto2011c,arpaia2018r,ArR:20}, and so on.

The ultimate goal of these approaches is to provide reliable and physically meaningful simulations for real-world applications,
while demanding minimal computational resources.
High-order methods are particularly suitable in this context, as they are able to achieve smaller errors within coarser discretizations.
Furthermore, an effective strategy for designing numerical methods with reduced errors is through structure-preserving techniques.
These techniques aim to replicate additional consistency conditions beyond those explicitly defined by the system of equations themselves.
For the SW equations, the focus is on preserving positive water heights, equilibrium or stationary states,
and implementing entropy conservation or dissipation methods.
First, the SW equations with source terms are known to admit a family of stationary solutions, which are characterized by a balance between flux divergence and source terms.
This concept is connected with that of a well-balanced (WB) discretization,
typically characterized by its ability to replicate one or more of these equilibria at the discrete level~\cite{berthon2016fully,berthon2022very,castro-pares-jsc20,GOMEZBUENO2021125820,gomez2021collocation,mantri2024fully,berberich2021high,chertock2018well}.
This WB property is crucial for complex, time-dependent simulations, as discretization errors due to the non-preservation of stationary regions could accumulate over time.
Second, in the context of flood simulations, it is a necessity to have provably positive discretizations, avoiding negative water heights.
To obtain a provably positive reconstruction in the context of high-order weighted essentially non-oscillatory~\cite{shu1998essentially}
(WENO) schemes, an effective positive limiter has been introduced and further developed in~\cite{zhang2010positivity,perthame1996positivity}. As proven in these references, this limiter achieves the preservation of positive reconstruction,
but it restricts the CFL condition, for classical SSPRK~\cite{gottlieb2001strong} schemes,
to the weight of the Gauss-Lobatto quadrature rule of the corresponding space accuracy (e.g., 1/12 for fifth-order schemes).
To circumvent this issue, unconditionally positivity preserving time-stepping strategies~\cite{meister2016positivity,ciallella2022arbitrary}
for the SW equations have been proposed,
based on a suitable reformulation of the finite volume semi-discretization in terms of production-destruction terms.
These approaches are based on the modified Patankar trick~\cite{patankar1980numerical,huang2019positivity,huang2018third,offner2020arbitrary}.
The linearly implicit nature of this approach allows for a relaxation of the aforementioned time-step constraint at a reasonable computational cost.

In this paper, we deal with the possibility of integrating the additional preservation of general static and moving equilibria
into the arbitrary high-order positivity preserving framework introduced in~\cite{ciallella2022arbitrary}.
To achieve this, we suitably modify the spatial discretization relying on ideas presented in~\cite{berthon2022very}.
In particular, we perform a convex blending between the original discretization and a WB one, able to exactly capture
general families of equilibria. The approach is able to tackle challenging flood simulations,
proving to be a good candidate for real-life applications.

The paper is structured as follows.
We first introduce the multidimensional SW system in \cref{se:SW_Equation}.
Then, the high-order WB positive numerical scheme is detailed in \cref{sec_Space},
where the space discretization is discussed, and in \cref{se_time_discretization}, where we present the production-destruction formulation in combination with high-order modified Patankar time schemes.
The results of the numerical validation are reported in \cref{se:numerics}.
Finally, \cref{se:summary} is left for conclusions and further developments.

\section{Shallow water equations}\label{se:SW_Equation}

%

%

The two-dimensional SW equations consist in a hyperbolic system of PDEs, extensively used in many applications to describe the behavior of water flows.
Their Eulerian formulation on a space domain $\Omega\subseteq\R^2$, assuming no friction and a time-independent bathymetry, reads
\begin{equation}\label{eq:CL}
\ddt\bu + \ddx\bF(\bu)+\ddy\bG(\bu) = \bS(x,y,\bu),\quad  \forall (x,y)\in \Omega, \quad\forall t \in [0,T_f],
\end{equation}
where conserved variables, fluxes and source term are respectively given by
\begin{equation}\label{eq:SWE}
	\begin{split}
&\bu=\begin{bmatrix} h \\ hu  \\ hv \end{bmatrix}\;,\;\;
\bS(x,y,\bu)=-gh\begin{bmatrix} 0 \\ \ddx b(x,y)  \\ \ddy b(x,y) \end{bmatrix},\\
&\bF(\bu)=\begin{bmatrix} hu\\
	hu^2+g\frac{h^2}{2}\\
	huv
\end{bmatrix},\;\;\;
\bG(\bu) =\begin{bmatrix}  hv\\huv \\hv^2+g\frac{h^2}{2} \end{bmatrix}\;,\;\;
\end{split}
\end{equation}
with $h$ being the water height,
$u$ and $v$ the velocity components of the flow along the $x$ and $y$ directions respectively,
$g$ the gravitational constant,
and $b(x,y)$ the bathymetry.
%
We also introduce the free surface water level $\eta\coloneqq h+b$, and the discharge variables along the two directions $x$ and $y$, defined as $q_x \coloneqq  h u$ and $q_y \coloneqq  h v$ respectively.

Notable properties of the SW equations, which have been drawing the interest of the scientific community in recent years and which play a central role in the context of this paper, are the positivity of the water height and the existence of non-trivial steady solutions.
In the context of numerical schemes preserving moving equilibria, one is interested in a detailed capturing of steady solutions satisfying
\begin{equation}\label{eq:Steady_State_General}
\ddt\bu\equiv 0 \Leftrightarrow \ddx\bF(\bu)+\ddy\bG(\bu) = \bS(x,y,\bu),\quad  \forall (x,y)\in \Omega\subseteq\R^2, \quad\forall t \in [0,T_f].
\end{equation}
The simplest and most known steady solution is the so-called ``lake at rest'' given by
\begin{equation}\label{eq:lake_at_rest}
u=v=0, \quad \eta \equiv\eta_0 \in \R^+_0, \quad \forall\, (x,y)\in \Omega,\, \forall t \in [0,T_f].
\end{equation}
Generally speaking, steady solutions are not known in closed-form and they are characterized by the analytical balance \eqref{eq:Steady_State_General}.
The smooth steady solutions tackled in this work are the pseudo-monodimensional states in the form
\begin{equation}
\left\{
\begin{aligned}
\frac{\partial}{\partial s} q_s = 0, \\
\frac{\partial}{\partial s} \left( \frac{q_s^2}{2h^2} + g(h+b) \right) = 0,
\end{aligned}
\right.
\label{eq:equilibriumODE}
\end{equation}
where $s$ is a general handle for the $x$ or the $y$ variable.
For more information on these steady solutions,
the reader is referred for instance to~\cite{michel2016well}.
%
For what follows, it is useful to define the so-called equilibrium variables
\begin{equation}
E_s(x,y,\mathbf{u}) = \begin{bmatrix}
q_s  \\
\frac{q_s^2}{2h^2} + g(h+b)
\end{bmatrix}.
\label{eq:equilibriumvariables}
\end{equation}
After \eqref{eq:equilibriumODE}, steady solutions are characterized by $E_s(x,y,\mathbf{u})$ being constant in space.

The system of PDEs under consideration is discretized using the Method of Lines (MOL), a numerical approach that treats space and time independently. In particular, space and time discretizations are the main focus of the next two sections.

\section{Well-balanced space discretization}\label{sec_Space}
This section is dedicated to the space discretization.
First, in \cref{sec:basic_space}, we describe our classical, non-well-balanced high-order discretization.
Then, \cref{subsec_Well_balance} is devoted to the generalization of a strategy
to achieve a high-order well-balanced (WB) discretization,
which was introduced in a one-dimensional setting in \cite{berthon2022very}.
Here, we generalize this technique for a two-dimensional WENO framework, applying the basic idea dimension by dimension.
The underlying principle consists in a simple blending between a high-order discretization and
a WB discretization to be used where a steady state is detected.
The main strengths of this approach are its low cost (no nonlinear equations need to be solved)
and its ease of use (it consists in multiplying the reconstruction by a suitable coefficient).
We emphasize that the resulting scheme will be able to capture and preserve
all the moving 1D steady solutions given by \eqref{eq:equilibriumODE},
and not just the so-called lake at rest solution, where velocity vanishes.
%

\subsection{Basic high-order discretization}\label{sec:basic_space}
The computational domain $\Omega$ is discretized in a Cartesian fashion via $N_x\times N_y$ non-overlapping control volumes
\begin{equation*}
\Omega_{i,j} = [\xin,\xip]\times[\yjn,\yjp],
\end{equation*}
with uniform spatial steps $\Delta x = \xip-\xin$ and $\Delta y = \yjp-\yjn$.

Finite volume methods are based on deriving a system of ordinary differential equations (ODEs) for the cell averages of the solution in each control volume $\Omega_{i,j}$
\begin{equation*}\label{eq:cell average}
\bU_{i,j}(t)\coloneqq \frac{1}{\Delta x \Delta y}\int_{\xin}^{\xip} \int_{\yjn}^{\yjp} \bu(x,y,t)\;\diff{x}\diff{y}.
\end{equation*}
The first step to obtain such a system is to integrate \eqref{eq:CL} over $\Omega_{i,j}$, thus getting
\begin{equation}\label{eq:evol cell average}
\ddo{\bU_{i,j}(t)}{t} + \frac{1}{\Delta x}(\bF_{\iip,j}(t)-\bF_{\iin,j}(t)) + \frac{1}{\Delta y}(\bG_{i,\jjp}(t)-\bG_{i,\jjn}(t)) = \bbS_{i,j}(t),
\end{equation}
where $\bbS_{i,j}$ is the source term average
\begin{equation*}\label{eq:source average}
\bbS_{i,j}(t)\coloneqq \frac{1}{\Delta x \Delta y}\int_{\xin}^{\xip} \int_{\yjn}^{\yjp} \bS(x,y,\bu)\;\diff{x}\diff{y},
\end{equation*}
and $\bF_{\iip,j}$ and $\bG_{i,\jjp}$ are the averages of the fluxes over the cell boundaries
%
\begin{align*}
\bF_{\iip,j}(t) &\coloneqq  \frac{1}{\Delta y}\int_{\yjn}^{\yjp}\bF(\bu(\xip,y,t))\;\diff{y},    \\
\bG_{i,\jjp}(t) &\coloneqq  \frac{1}{\Delta x}\int_{\xin}^{\xip}\bG(\bu(x,\yjp,t))\;\diff{x}  .
\end{align*}
%
So far, Equation \eqref{eq:evol cell average} has been exactly derived from \eqref{eq:CL}.
However, in order to obtain the numerical scheme, we need to discretize the fluxes and the source averages.

To that end, we rely on the following ingredients:
a high-order reconstruction of the conservative variables in each control volume (WENO \cite{shu1998essentially} in our case), consistent quadrature formulas to discretize all integrals (Gauss-Legendre with $Q$ points in our case), and suitable numerical fluxes to compute the fluxes in the boundary integrals (described later on).
In the remainder of this section, we drop the time dependency to shorten notation.


Let us first focus on the discretization of the fluxes averages and, more in detail, on $\bF_{\iip,j}$, as $\bG_{i,\jjp}$ is obtained similarly.
Once the reconstruction in each control volume has been performed, at each quadrature point $y_q \in [y_{j-1/2},y_{j+1/2}]$ of each edge $\xip$ we have two high-order reconstructed values for $\bu$, corresponding to $\xip^L$ and $\xip^R$,
which will be referred to as the left and right high-order extrapolated values
\begin{equation*}
\bu^L_{\iip,q} = \bu^{\text{HO}}(\xip^L,y_q)
\text{\qquad and \qquad}
\bu^R_{\iip,q} = \bu^{\text{HO}}(\xip^R,y_q).
\end{equation*}
By applying a consistent quadrature rule, the flux in the $x$-direction reads
\begin{equation*}\label{eq:flux with GP}
\bF_{\iip,j} \approx \sum_{q=1}^{Q} w_q \hat\bF(\bu_{\iip,q}^L,\bu_{\iip,q}^R),
\end{equation*}
where $\hat\bF$ is a consistent numerical flux, and $w_q$ is the normalized quadrature weight associated to the quadrature node $y_q$.
The choice of $\hat\bF$ is discussed in \cref{subsec_Well_balance}.

The high-order source term averages are computed as
\begin{equation*}
\bbS_{i,j} \approx \sum_{q=1}^{Q}\sum_{p=1}^{Q} w_q w_p \bS(x_q,y_p,\bu^{HO}(x_q,y_p)),
\end{equation*}
with a surface quadrature obtained as the tensor product of the classical 1D quadrature used for the edges and $\bu^{HO}$ being the local reconstruction of the solution in the cell.
%
%

Despite its robustness in capturing discontinuities, while minimizing the oscillations, the WENO reconstruction may provide some negative reconstructed values for the water height, especially close to dry regions.
Such negative water heights are not physically admissible,
and in fact will immediately lead to the simulation crashing.
In order to avoid such an issue, we adopt for the water height reconstruction the
positivity limiter introduced in~\cite{perthame1996positivity} and further discussed in~\cite{zhang2010positivity}.

As shown in \cite{xing2014survey}, provable positivity preservation for the water height, in the context of this framework, is subjected to severe CFL constraints, when adopting standard time integration techniques.
In particular, assuming a simple forward Euler time-stepping and a Lax-Friedrichs numerical flux, the limit CFL guaranteeing positivity preservation is $\text{CFL}^{\text{FE}}\coloneqq w^{\text{Lobatto}}_{1}$, where $w^{\text{Lobatto}}_{1}$ is the first weight of the adopted high-order Gauss-Lobatto quadrature rule.
This corresponds to $\text{CFL}^{\text{FE}}=1/12$ for a quadrature of order~$5$.
The restriction gets even worse as the order of accuracy increases, e.g., we have $\text{CFL}^{\text{FE}}=1/20$ for a quadrature of order~$7$.
The adoption of high-order SSPRK methods slightly relaxes the constraint, but not significantly.
Indeed, for instance, using the SSPRK$(5,4)$ discretization relaxes the condition to $\text{CFL}^{\text{SSPRK}(5,4)} = 1.508 \,\text{CFL}^{\text{FE}}$.
The adopted time discretization, described in \cref{se_time_discretization}, allows us to drop such limitations and to run simulations at any CFL without violating the positivity constraint on the water height.
Due to the explicit nature of the time scheme used for the discharge equations, however, the (far less restrictive) stability constraint $\text{CFL}\leq 1$ of explicit schemes applies.


\subsection{Well-balanced blending}
\label{subsec_Well_balance}

We now describe the WB strategy, which makes possible the capture of steady states characterized by constant equilibrium variables~\eqref{eq:equilibriumvariables}. 
The key idea comes from the following remark: for the simulation of a steady solution, a well-balanced scheme is exact, and therefore has a better accuracy than any high-order scheme.
For unsteady simulations, high-order schemes are more accurate, and should be used whenever the solution is not steady.
To achieve a seamless switch between high-order and well-balanced schemes, we propose a simple blending between the two.
This blending is performed according to a suitable steady solution indicator, defined below.

For simplicity, we only derive the reconstruction along the $x$-direction.
The extension to the $y$-direction is easily performed following a dimension by dimension approach.
We replace the reconstructed variables at the interfaces by the convex combination between the high-order extrapolated values and the cell averages
\begin{equation}
    \label{eq:modified_reconstruction}
    \begin{aligned}
\tilde{\mathbf{u}}^L_{\iip,q} &= (1 - \theta_{\iip,j}) \mathbf{U}_{i,j} + \theta_{\iip,j} \mathbf{u}^L_{\iip,q}, \\
\tilde{\mathbf{u}}^R_{\iip,q} &= (1 - \theta_{\iip,j}) \mathbf{U}_{i+1,j} + \theta_{\iip,j} \mathbf{u}^R_{\iip,q},
    \end{aligned}
\end{equation}
where $\theta_{\iip,j}$ is a steady state indicator.
On the one hand, it should vanish when the equilibrium variables~\eqref{eq:equilibriumvariables} are constant in space; in this case, the modified reconstructed values $\tilde{\mathbf{u}}^L_{\iip,q}$ are equal to the cell averages $\mathbf{U}_{i,j}$.
On the other hand, when far from any equilibrium, $\tilde{\mathbf{u}}^L_{\iip,q}$ should be an approximation of order $P$, where $P$ is the order of the discretization (herein, $P=5$).

Following~\cite{berthon2022very}, we define $\theta_{\iip,j}$ by
\begin{equation*}
\theta_{\iip,j} = \frac{\varepsilon_{\iip,j}}{\varepsilon_{\iip,j} + \left(\frac{\Delta x}{C_{\iip,j}}\right)^P},
\end{equation*}
with
\begin{equation*}
\varepsilon_{\iip,j} \coloneqq  \| E_x(x_{i+1},y_j,\theta_{\iip,j}) - E_x(x_{i},y_j,\theta_{\iip,j}) \|,
\end{equation*}
where $C_{\iip,j}$ is a quantity independent of $\Delta x$, which is here chosen, at a given time iteration, as the time residual difference at the previous iteration
\begin{equation*}
C_{\iip,j} \coloneqq  \frac12\left( \frac{\bU^{n}_{i+1,j} - \bU^{n-1}_{i+1,j}}{\Delta t} + \frac{\bU^{n}_{i,j} - \bU_{i,j}^{n-1}}{\Delta t} \right).
\end{equation*}
We remark that, at equilibrium $C_{\iip,j}\to 0$, hence $\theta_{\iip,j}\to 0$ as well resulting in the low-order WB reconstruction.

Similarly, the source term discretization is defined as
\begin{equation*}
    \begin{aligned}
        \tilde{\bbS}_{i,j}
        & = \frac12\left( \theta_{\iin,j} + \theta_{\iip,j} \right) \bbS_{i,j} \\
        & + \frac12 \left(  \left(1- \theta_{\iin,j}\right)\bbS^{\text{WB}}_{\iin,j}  + \left(1-\theta_{\iip,j}\right)\bbS^{\text{WB}}_{\iip,j} \right),
    \end{aligned}
\end{equation*}
where $\bbS^{\text{WB}}_{\iin,j}$ and $\bbS^{\text{WB}}_{\iip,j}$ represent WB discretizations of the source term at the interfaces described in~\cite{michel2016well}.

In order to ensure stability in the context of unsteady wet-dry simulations, we found experimentally useful to introduce, with respect to the classical approach~\cite{berthon2022very,michel2016well}, a similar convex combination in the flux definition
\begin{equation*}
    \begin{aligned}
        \hat{\bF}(\mathbf{u}^L_{\iip,q},\mathbf{u}^R_{\iip,q})
        &=  (1 - \theta_{\iip,j})\hat{\bF}^{\text{WB}}(\mathbf{u}^L_{\iip,q},\mathbf{u}^R_{\iip,q}) \\
        &+ \theta_{\iip,j} \hat{\bF}^{\text{LF}}(\mathbf{u}^L_{\iip,q},\mathbf{u}^R_{\iip,q}),
    \end{aligned}
\end{equation*}
where $\hat{\bF}^{\text{WB}}$ represents the WB approximate Riemann solver presented in~\cite{michel2016well}, while $\hat{\bF}^{\text{LF}}$ is a robust local Lax-Friedrichs numerical flux reading
\begin{equation*}
\hat{\bF}^{\text{LF}}(\mathbf{u}^L,\mathbf{u}^R) = \frac12\left( \bF(\mathbf{u}^R) +\bF(\mathbf{u}^L)\right) - \frac12 s_{max} \left( \mathbf{u}^R + \mathbf{u}^L \right),
\end{equation*}
where $s_{max}$  is the spectral radius of the normal flux Jacobian of system~\eqref{eq:CL}.

\begin{remark}\label{rmk:WBLO}
The discretized terms $\bbS^{\text{WB}}_{\iin,j}$, $\bbS^{\text{WB}}_{\iip,j}$ and $\hat{\bF}^{\text{WB}}$ are designed in such a way to guarantee an exact equilibrium with respect to steady states in the form~\eqref{eq:equilibriumvariables}, when taking in input the cell averages.
The reader can easily verify that, when a steady state of this type is considered, then all $\theta_{\iip,j}$ are equal to 0 and the scheme reduces to the WB version.
Indeed, the modified reconstruction \eqref{eq:modified_reconstruction} degenerates to the cell averages.
This means that, despite guaranteeing an exact capturing of the steady states, the basic WB discretization is directly based on cell averages without any reconstruction, and it is, therefore, only first order accurate in general \cite{berthon2022very}.
\end{remark}

\begin{remark}
    We emphasize an important property of the proposed strategy: the steady solution indicator is defined such that the nonlinear system \eqref{eq:equilibriumODE} never has to be solved.
    Instead, it merely relies on evaluating the equilibrium variables \eqref{eq:equilibriumvariables} at the cell interfaces.
\end{remark}

\begin{remark}
It should be noticed that the proposed scheme is well-balanced when $\theta_{\iip,j}$ goes to $0$. Numerically speaking, this consists in defining a low enough threshold ($10^{-10}$ in the numerical experiments) to set $\theta_{\iip,j}= 0$.
\end{remark}

\section{Unconditionally positive time discretization}\label{se_time_discretization}

In this section, we describe the time-stepping strategy, which consists in a slight modification of arbitrary high-order deferred correction (DeC) methods for ODEs~\cite{abgrall2017high}.
In particular, the water height update is reinterpreted as a Production-Destruction System (PDS) and then the modified Patankar trick is applied in order to achieve unconditional preservation of its positivity as in \cite{ciallella2022arbitrary}.
Both DeC methods and Patankar trick have a long history.
In particular, for more information on DeC the interested reader is referred to~\cite{dutt2000dec,micalizzi2023new,micalizzi2023efficient,veiga2024improving,veiga2021dec},
while Patankar (and modified Patankar) tricks are detailed in~\cite{patankar1980numerical,offner2020arbitrary,burchard2003high,burchard2005application,izgin2022stability,kopecz2018order,kopecz2018unconditionally,huang2019positivity,huang2018third,meister2014unconditionally,ortleb2017patankar}.

\subsection{Deferred Correction method}
\label{sec:DeC}

To introduce the DeC method, let us consider the Cauchy problem
\begin{equation}\label{eq:initial_prob}
\begin{cases}
\frac{d}{dt}\uvec{c}(t) = \uvec{H}(t,\uvec{c}(t)), \quad t\in[0,T_f],\\
\uvec{c}(0) = \uvec{c}_0,
\end{cases}
\end{equation}
where $\uvec{c}:[0,T_f]\to \R^{N_c}$ is the unknown solution, with $N_c$ components, and $\uvec{H}:[0,T_f] \times \R^{N_c} \to \R^{N_c}$ is a given function satisfying the classical smoothness assumptions, which guarantee the existence of a unique solution to the Cauchy problem \eqref{eq:initial_prob}.
As is customary in the context of one-step methods, we focus on a generic interval $[t_n,t_{n+1}]$ of size $\Delta t\coloneqq t_{n+1}-t_n$ and, given $\uvec{c}_n\approx\uvec{c}(t_n)$, we seek an approximation $\uvec{c}_{n+1}$ of~$\uvec{c}(t_{n+1})$.

Following~\cite{abgrall2017high,micalizzi2023new}, we introduce $M+1$ subtimenodes $t^m$ in the interval $[t_n,t_{n+1}]$, which are such that
\begin{equation*}
t_n=t^0<t^1<\dots<t^M=t_{n+1}.
\end{equation*}
The DeC method under consideration consists in an explicit fixed point iterative procedure
to compute the approximation of $c$ at all subtimenodes simultaneously.
The update formula is given by
\begin{equation}
\uvec{c}^{m,(p)}\coloneqq \uvec{c}^0+\Delta t\sum_{\ell=0}^M \theta_\ell^m \uvec{H}(t^{\ell},\uvec{c}^{\ell,(p-1)}),\quad m=1,\dots,M,\quad p\geq 1,
\label{eq:DeC}
\end{equation}
where $\uvec{c}^{m,(p)}$ is the approximation of the solution in the subtimenode $t^m$ obtained at the $p^\text{th}$ iteration and, for each $m$, the coefficients $(\theta_\ell^m)_{\ell \in \{0,\dots,M\}}$ are the normalized weights of the high-order quadrature formula over $[t^0,t^m]$ associated to the subtimenodes.
In particular, in the previous update formula, we set $\uvec{c}^{m,(p)}=\uvec{c}^0\coloneqq \uvec{c}_n$ whenever $m=0$ or $p=0$.
One can show that, for small enough $\dt$, the iterative process converges. Furthermore, the order of accuracy of $\uvec{c}^{M,(p)}$ with respect to $\uvec{c}(t_{n+1})$ is $\min{(p,R)}$, i.e., each iteration corresponds to an increase in the order of accuracy by one, until a saturation value $R$, which depends on the number and on the distribution of the adopted subtimenodes.
For example, evenly spaced subtimenodes lead to $R=M+1$, while Gauss-Lobatto subtimenodes yield $R=2M$.
In this paper, we use Gauss-Lobatto subtimenodes.
Therefore, the optimal way of reaching order $P$ is to perform $P$ fixed-point iterations with $M+1$ subtimenodes, where $M=\lceil\frac{P}{2}\rceil$.

Hence, the arbitrarily high-order time integration method presented in this section, combined with the space discretization described in \cref{sec_Space}, defines an arbitrarily high-order, fully well-balanced framework for the numerical solution of the SW equations~\eqref{eq:CL} -- \eqref{eq:SWE}.
However, at this level, nothing can be said, in general, about the positivity of the water height.
In the next subsection, we present the modification to be performed in the time integration of the water height, guaranteeing unconditional positivity.

\subsection{Modified Patankar DeC method}

In this section, we first focus, in \cref{sec:PDS}, on the unconditionally positive time integration of a specific class of ODEs, namely Production-Destruction Systems (PDSs).
Then, we describe in \cref{sec:SW_as_PDS} how to apply these notions to the SW equations.

\subsubsection{Unconditionally positive time integration of PDSs}
\label{sec:PDS}
PDSs are systems of ODEs characterized by the following structure
\begin{equation*}
	\begin{cases}
		\dfrac{d}{dt}c_\alpha=\sum_{\beta=1}^{N_c}p_{\alpha,\beta}(\uvec{c}) -\sum_{\beta=1}^{N_c}d_{\alpha,\beta}(\uvec{c}), \quad \alpha=1,\dots,N_c,\\
		\uvec{c}(0)=\uvec{c}_0, \vphantom{\dfrac{d}{dt}}
	\end{cases}
\end{equation*}
where $\uvec{c} = (c_\alpha)_{\alpha \in \{1,\dots,N_c\}}$,
and where $p_{\alpha,\beta}$ and $d_{\alpha,\beta}$ are real non-negative Lipschitz-continuous functions from $\R^{N_c}$ to $\R_{+}^*$.

More specifically, we are interested in a subfamily of PDSs fulfilling two extra constraints: conservation and positivity.
A PDS is said to be conservative if, $\forall \alpha, \beta \in \{1,\dots,N_c\}$ and $\forall\uvec{c}\in \R^{N_c}$, we have
$\Prod_{\alpha,\beta}(\uvec{c}) = \dest_{\beta,\alpha}(\uvec{c})$, thus implying
\begin{equation}
	\label{eq:PDS_conservation}
	\sum_{\alpha=1}^{N_c}c_\alpha(t)=\sum_{\alpha=1}^{N_c}c_\alpha(0), \quad \forall t \in [0,T_f].
\end{equation}
A PDS is said to be positive if, starting by a positive initial condition, we get a positive evolution of all the components, i.e.,
\begin{equation}
	\label{eq:PDS_positivity}
	\uvec{c}(0)>0 \implies \uvec{c}(t)>0, \quad \forall t\in[0,T_f],
\end{equation}
where the comparison operator, applied to vectors, is meant to be applied to each scalar component.

Conservative and positive PDSs arise in many applications and many numerical methods have been developed to preserve such properties.
A successful approach, in this context, is given by the (modified) Patankar trick~\cite{patankar1980numerical,burchard2003high}, which is based on the introduction of some weights on the production and destruction terms.
In particular, the application of the modified Patankar trick to the DeC scheme (mPDeC) \cite{offner2020arbitrary}
is characterized by replacing \eqref{eq:DeC} with the following update
\begin{equation*}
	\label{eq:mPDeC}
	c_\alpha^{m,(p)}=c_\alpha^0 +\Delta t \sum_{\ell=0}^M \theta_\ell^m \left(   \sum_{\beta=1}^{N_c}  \Prod_{\alpha,\beta}(\uvec{c}^{\ell,(p-1)}) \frac{c^{m,(p)}_{\gamma(\beta,\alpha, \theta_\ell^m)}}{c_{\gamma(\beta,\alpha, \theta_\ell^m)}^{m,(p-1)}} - \sum_{\beta=1}^{N_c} \dest_{\alpha,\beta}(\uvec{c}^{\ell,(p-1)})  \frac{c^{m,(p)}_{\gamma(\alpha,\beta, \theta_\ell^m)}}{c_{\gamma(\alpha,\beta, \theta_\ell^m)}^{m,(p-1)}} \right),
\end{equation*}
where $\uvec{c}^{m,(p)}=\uvec{c}^{0}\coloneqq \uvec{c}_n$ whenever $m=0$ or $p=0$, and $\gamma$ is a switch function defined as
\begin{equation*}
    \gamma(\alpha,\beta,\theta)\coloneqq \begin{cases}
        \alpha, & \text{if }\theta \geq 0,\\
        \beta, & \text{if }\theta <0.
    \end{cases}
\end{equation*}
For guidelines concerning the number $P$ of iterations to be performed, and the associated accuracy, the reader is referred to the discussion regarding the standard DeC scheme at the end of \cref{sec:DeC}.
The mPDeC method is positive and conservative, i.e., it satisfies
\begin{equation*}
	\sum_{\alpha=1}^{N_c}c_{\alpha,n+1}=\sum_{\alpha=1}^{N_c}c_{\alpha,n}
    \text{\qquad and \qquad}
	\uvec{c}_n>0 \implies \uvec{c}_{n+1}>0,
\end{equation*}
which are nothing but natural translations, at the discrete level, of the continuous constraints \eqref{eq:PDS_conservation} and~\eqref{eq:PDS_positivity}.

Moreover, the method is linearly implicit and can be recast in compact form as
\begin{equation*}
	\label{eq:mPDeC_system}
	\mathbb M \uvec{c}^{m,(p)} = \uvec{c}_n,
\end{equation*}
where the matrix $\mathbb M$ is  defined as
{\small	\begin{equation}\label{eq:matrixMPDeC}%
		\begin{split}
		&\mathbb M(\underline{\uvec{c}}^{(p-1)},m)_{\alpha,\beta}= \\
		&\begin{cases}
			1+\Delta t \sum\limits_{\ell=0}^M \sum\limits_{\substack{k=1 \\ k\neq \alpha}}^{N_c}  \frac{\theta_\ell^m}{c_\alpha^{m,(p-1)}}   \left(  d_{\alpha,k}(\uvec{c}^{\ell,(p-1)})  \chi_{\lbrace \theta^m_\ell\geq 0\rbrace}  - p_{i,k} (\uvec{c}^{\ell,(p-1)})\chi_{\lbrace \theta^m_\ell<0 \rbrace} \right),  \,&\text{for } \beta=\alpha,\\
			-\Delta t \sum\limits_{\ell=0}^M \frac{\theta_l^m}{c_\beta^{m,(p-1)}} \left(  p_{\alpha,\beta}(\uvec{c}^{\ell,(p-1)})\chi_{\lbrace \theta^m_\ell \geq 0\rbrace} - d_{\alpha,\beta} (\uvec{c}^{\ell,(p-1)})\chi_{\lbrace \theta^m_\ell<0\rbrace} \right),  \, &\text{for } \beta\ne \alpha,
		\end{cases}
		\end{split}
	\end{equation}}%
with $\chi_{\lbrace \cdot \rbrace}$ the indicator function, i.e., a switch with value equal to $1$ if the argument condition is true, $0$ otherwise.
One can prove that the matrix is column diagonally dominant and hence invertible. Furthermore, it is possible to show that, for any $\uvec{b}>0$, the solution to $\mathbb M \uvec{c}=\uvec{b}$ is such that $\uvec{c}>0.$
At the implementation level, the system is solved though the Jacobi method, which is provably convergent due to the fact that~$\mathbb M$ is column diagonally dominant.
Moreover, in order to avoid divisions by zero, as in~\cite{meister2016positivity,offner2020arbitrary}, the following mollification of the ratios in the matrix \eqref{eq:matrixMPDeC} is considered
\begin{equation*}
	\frac{n}{d}\approx \begin{cases}
		0, & \text{ if } d < 10^{-8},\\
		\frac{2d\cdot n}{d^2+\max\lbrace d^2,10^{-8} \rbrace}, & \text{ if } d \geq 10^{-8}.
	\end{cases}
\end{equation*}

Further details are omitted to avoid lengthening the paper.
However, they are thoroughly discussed in \cite{ciallella2022arbitrary}.

We now explain how the presented notions can be applied to the finite volume semi-discretization of the SW equations.

\subsubsection{Application to the Shallow Water equations}
\label{sec:SW_as_PDS}

The key idea is to reinterpret the water height semi-discretization as a PDS and to apply the modified Patankar trick to the water height DeC update, while performing a standard DeC time-stepping on the updates of the discharge in the $x$- and $y$-directions.
From \eqref{eq:evol cell average}, we note that each cell communicates with the neighboring cells, sharing common edges, via numerical fluxes.
Thus, in such a context, the components $c_\alpha$ are given by the water height averages $h_{i,j}$, with indices $\alpha$ identified as couples $[i,j]$, while the production and destruction terms are given by the associated numerical fluxes. Let us recall that the water height equation has no source term contribution.
%
%

\begin{figure}[!ht]
	\centering
	\scalebox{1.1}{
		\begin{tikzpicture}

			\draw (0.,2) -- (2,2) -- (2,3.5) -- (0.,3.5) -- (0.,2);
			\draw (2,2) -- (3.5,2) -- (3.5,3.5) -- (2,3.5) -- (2,2);
			\draw (3.5,2) -- (5.5,2) -- (5.5,3.5) -- (3.5,3.5) -- (3.5,2);
			\draw (2,0.5) -- (3.5,0.5) -- (3.5,2) -- (2,2) -- (2,0.5);
			\draw (2,3.5) -- (3.5,3.5) -- (3.5,5) -- (2,5) -- (2,3.5);

			\node [black,scale=0.8] at (0.9,2.75) {$\Omega_{i-1,j}$};
			\node [black,scale=0.8] at (2.75,2.75) {$\Omega_{i,j}$};
			\node [black,scale=0.8] at (4.6,2.75) {$\Omega_{i+1,j}$};
			\node [black,scale=0.8] at (2.75,4.25) {$\Omega_{i,j+1}$};
			\node [black,scale=0.8] at (2.75,1.25) {$\Omega_{i,j-1}$};

			\node [red, scale=0.7, anchor=south] at (1.9,2.7) {$\hat{\bF}^{(1)}_{i-1/2,j}(\bU)$};

			\draw [-stealth,red] (1.55,2.7) -- (2.25,2.7);

			\node [red, scale=0.7, anchor=south] at (3.6,2.7) {$\hat{\bF}^{(1)}_{i+1/2,j}(\bU)$};

			\draw [stealth-,red] (3.95,2.7) -- (3.25,2.7);

			\node [red, scale=0.7, anchor=base west] at (2.85,2.05) {$\hat{\bG}^{(1)}_{i,j-1/2}(\bU)$};

			\draw [-stealth,red] (2.75,1.65) -- (2.75,2.35);

			\node [red, scale=0.7, anchor=base west] at (2.85,3.65) {$\hat{\bG}^{(1)}_{i,j+1/2}(\bU)$};

			\draw [-stealth,red] (2.75,3.15) -- (2.75,3.85);

	\end{tikzpicture}}\hfill
	\scalebox{1.1}{
		\begin{tikzpicture}

			\draw (0.,2) -- (2,2) -- (2,3.5) -- (0.0,3.5) -- (0.,2);
			\draw (2,2) -- (3.5,2) -- (3.5,3.5) -- (2,3.5) -- (2,2);
			\draw (3.5,2) -- (5.5,2) -- (5.5,3.5) -- (3.5,3.5) -- (3.5,2);
			\draw (2,0.5) -- (3.5,0.5) -- (3.5,2) -- (2,2) -- (2,0.5);
			\draw (2,3.5) -- (3.5,3.5) -- (3.5,5) -- (2,5) -- (2,3.5);

			\node [black,scale=0.8] at (0.9,2.75) {$\Omega_{i-1,j}$};
			\node [black,scale=0.8] at (2.75,2.75) {$\Omega_{i,j}$};
			\node [black,scale=0.8] at (4.6,2.75) {$\Omega_{i+1,j}$};
			\node [black,scale=0.8] at (2.75,4.25) {$\Omega_{i,j+1}$};
			\node [black,scale=0.8] at (2.75,1.25) {$\Omega_{i,j-1}$};

			\node [red, scale=0.7, anchor=base] at (1.9,2.48) {$p_{[i,j],[i-1,j]}$};
			\node [red, scale=0.7, anchor=south] at (1.9,2.82) {$d_{[i,j],[i-1,j]}$};

			\draw [-stealth,red] (1.55,2.7) -- (2.25,2.7);
			\draw [stealth-,red] (1.55,2.8) -- (2.25,2.8);

			\node [red, scale=0.7, anchor=base] at (3.6,2.48) {$p_{[i,j],[i+1,j]}$};
			\node [red, scale=0.7, anchor=south] at (3.6,2.82) {$d_{[i,j],[i+1,j]}$};

			\draw [-stealth,red] (3.95,2.7) -- (3.25,2.7);
			\draw [stealth-,red] (3.95,2.8) -- (3.25,2.8);

			\node [red, scale=0.7, anchor=base east, xshift=-1] at (2.65,1.75) {$p_{[i,j],[i,j-1]}$};
			\node [red, scale=0.7, anchor=base west] at (2.85,1.75) {$d_{[i,j],[i,j-1]}$};

			\draw [-stealth,red] (2.7,1.65) -- (2.7,2.35);
			\draw [stealth-,red] (2.8,1.65) -- (2.8,2.35);

			\node [red, scale=0.7, anchor=base east, xshift=-1] at (2.65,3.7) {$d_{[i,j],[i,j+1]}$};
			\node [red, scale=0.7, anchor=base west] at (2.85,3.7) {$p_{[i,j],[i,j+1]}$};

			\draw [-stealth,red] (2.7,3.15) -- (2.7,3.85);
			\draw [stealth-,red] (2.8,3.15) -- (2.8,3.85);

	\end{tikzpicture}}
	\caption{Sketch of the PDS structure for the control volume $\Omega_{i,j}$.}\label{fig:omega_ij}
\end{figure}
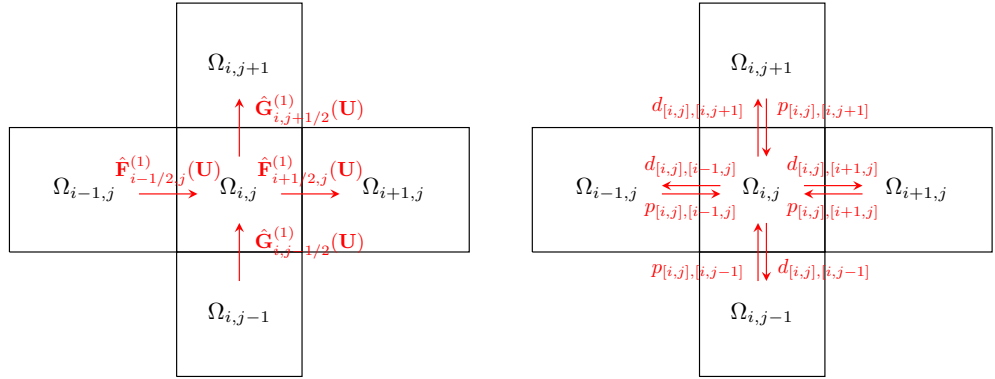

Considering all the neighbors to the cell $[i,j]$ (i.e., cells sharing an edge with cell $[i,j]$), one can define the following production and destruction terms
\begin{equation}\label{eq:PDSforWENO}
	\begin{split}
		p_{[i,j],[i-1,j]}({\bU})= + \frac{1}{\Delta x}  \hat{\bF}^{(1)}_{i-1/2,j} (\bU)^+, \quad d_{[i,j],[i-1,j]}({\bU})=- \frac{1}{\Delta x} \hat{\bF}^{(1)}_{i-1/2,j} (\bU)^-,\\
		p_{[i,j],[i+1,j]}({\bU})= - \frac{1}{\Delta x}  \hat{\bF}^{(1)}_{\iip,j}  (\bU)^-, \quad d_{[i,j],[i+1,j]}({\bU})=+ \frac{1}{\Delta x} \hat{\bF}^{(1)}_{\iip,j}  (\bU)^+,\\
		p_{[i,j],[i,j-1]}({\bU})= + \frac{1}{\Delta y}  \hat{\bG}^{(1)}_{i,j-1/2} (\bU)^+, \quad d_{[i,j],[i,j-1]}({\bU})=- \frac{1}{\Delta y} \hat{\bG}^{(1)}_{i,j-1/2}  (\bU)^-,\\
		p_{[i,j],[i,j+1]}({\bU})= - \frac{1}{\Delta y}  \hat{\bG}^{(1)}_{i,j+1/2} (\bU)^-, \quad d_{[i,j],[i,j+1]}({\bU})=+ \frac{1}{\Delta y} \hat{\bG}^{(1)}_{i,j+1/2}  (\bU)^+,
	\end{split}
\end{equation}
where the superscripts $\,^+$ and $\,^-$ respectively represent the positive and the negative part, while the superscript~$^{(1)}$ represents the first component of the numerical fluxes.
%
%
These production and destruction terms, as well as their relationships \eqref{eq:PDSforWENO} with the numerical fluxes, are sketched in \cref{fig:omega_ij}.

In light of the previous discussion, it is therefore easy to apply the modified Patankar trick to the DeC update of the water height.
Let us remark that the strategy provides unconditional positivity of the water height with respect to the time step $\dt$.
This results in great computational advantages with respect to standard explicit time integration techniques subjected to the typical positivity-preserving CFL constraints.
%
Further details, including a detailed description of a possible implementation, can be found in~\cite{ciallella2022arbitrary}.


%

\section{Numerical results}\label{se:numerics}

In this section, we report the results of several numerical experiments demonstrating the good properties of the scheme, including its robustness.
In particular, the tests are meant to verify the high-order accuracy in \cref{sec:vortex}, the WB property for both stationary and moving equilibria (in \cref{sec:lake_at_rest,sec:simu_moving}), and the ability to deal with tough flood simulations involving dry areas in \cref{sec:simu_flooding}.
We assume $g=9.81$ unless otherwise specified.
Let us remark that the basic ingredients of the scheme allow us to reach arbitrarily high-orders of accuracy. Here, we focus on the fifth order version.
%

\subsection{Unsteady vortex}
\label{sec:vortex}


Through this test \cite{ricchiuto2021analytical}, we verify the high-order accuracy of the space and time discretizations, without considering the source term for the moment. To that end, we set $b\equiv 0$.
Therefore, this test is meant to verify the high-order accuracy of the flux discretization; the high-order accuracy of the source term discretization will be checked in a later test.

The considered computational domain is the square $\Omega\coloneqq [0,3]\times[0,3]$, and the vortex is given by a perturbation $\delta$ of a homogeneous background field $(h_0,u_0,v_0)\coloneqq (1,2,3)$.
Let us define the variable $r(x,y,t)\coloneqq \sqrt{(x-x_c(t))^2+(y-y_c(t))^2}$, expressing the distance between $(x,y)$ and vortex center $(x_c(t),y_c(t))\coloneqq (1.5,1.5)+(u_0 t, v_0 t)$.

The water height is then given by $h(r)\coloneqq h_0 + \delta h(r)$, with
\begin{equation*}\label{eq:movingvortexh}
	\delta h(r) \coloneqq  - \gamma
	\begin{cases}
		\exp \left( -\dfrac{1}{\arctan^3(1-r^2)} \right), & \text{ if } r < 1,                       \\
		0,                                                & \text{ otherwise}, \vphantom{\dfrac 1 2}
	\end{cases}
\end{equation*}
where $\gamma\coloneqq 0.1$ is the vortex amplitude.
The velocity field, defined by $(u,v)\coloneqq (u_0,v_0)+ (\delta u,\delta v)$,
is characterized by the following perturbation
\begin{equation*}\label{eq:movingvortexu}
	\begin{pmatrix}
		\delta u \\ \delta v
	\end{pmatrix}
	=
	\sqrt{\frac{g}{r}\,\frac{\partial h}{\partial r} }
	\begin{pmatrix}
		y-y_c \\ -(x-x_c)
	\end{pmatrix},
\end{equation*}
where $\frac{\partial h}{\partial r}$ is the derivative of $h$ with respect to $r$,
which depends only on the radial distance from the center of the vortex
%
%
\begin{equation*}
	\frac{\partial h}{\partial r}(r) =
	\begin{cases}
		\dfrac{6\,\gamma\, r\,\exp(-\frac{1}{\arctan^3(1-r^2)})}{\arctan^4(r^2-1)(1 + (r^2-1)^2)}, & \text{ if } r < 1,                       \\
		0,                                                                                      & \text{ otherwise}. \vphantom{\dfrac 1 2}
	\end{cases}
\end{equation*}
We assume periodic boundary conditions and a final time $T_f\coloneqq 0.1$.
It is important to highlight the fact that this solution is $\mathcal{C}^\infty$, which is a fundamental property for testing arbitrarily high-order schemes \cite{ricchiuto2021analytical}.

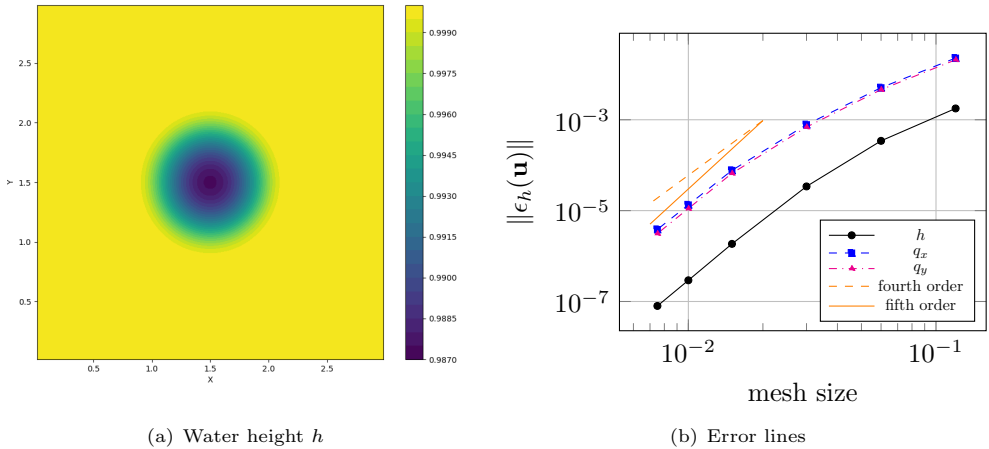
\begin{figure}
	\centering
	\subfigure[Water height $h$]{
		\includegraphics[width=0.47\textwidth,trim={1.5cm 0cm 2.5cm 2.3cm},clip]{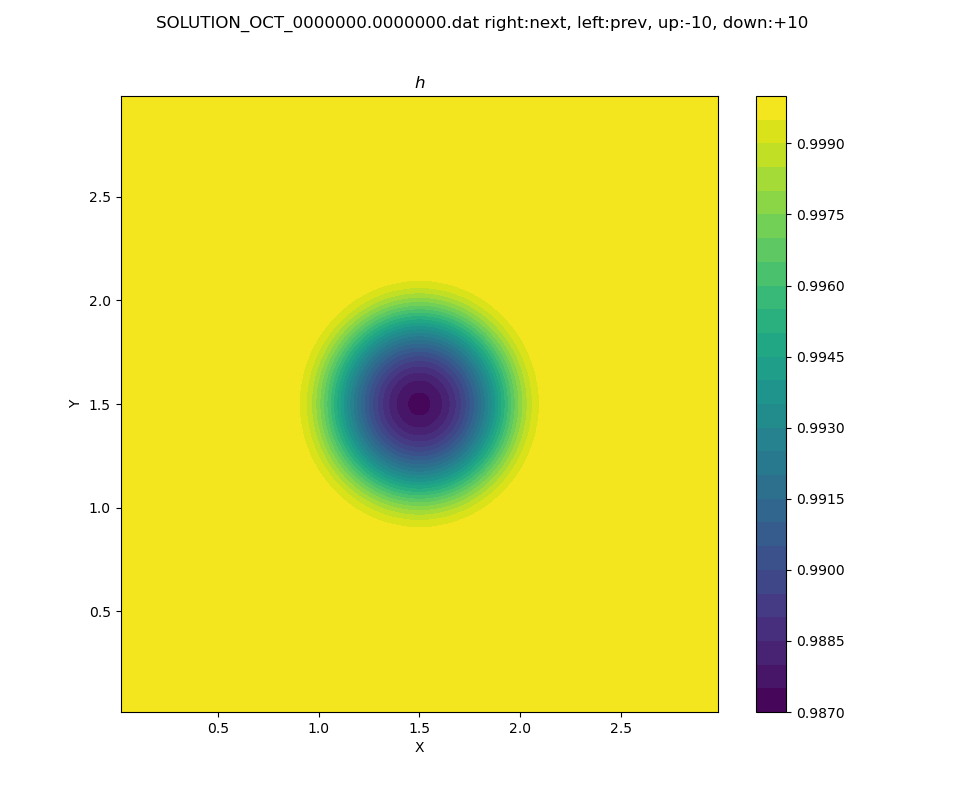}
	}\hfill
	\subfigure[Error lines]{
		\begin{tikzpicture}
			\begin{axis}[
					xmode=log, ymode=log,
					grid=major,
					xlabel={mesh size},
					ylabel={$\|\epsilon_h(\mathbf{u})\|$},
					xlabel shift = 1 pt,
					ylabel shift = 1 pt,
					legend pos= south east,
					legend style={nodes={scale=0.6, transform shape}},
					width=.49\textwidth
				]
				\addplot[mark=otimes*,mark size=1.3pt,black] table [y=h, x=size]{vortex_mPDec5_with_orders.dat};
				\addlegendentry{$h$}

				\addplot[mark=square*,mark size=1.3pt,blue,dashed] table [y=hu, x=size]{vortex_mPDec5_with_orders.dat};
				\addlegendentry{$q_x$}

				\addplot[mark=triangle*,mark size=1.3pt,magenta,dashdotted] table [y=hv, x=size]{vortex_mPDec5_with_orders.dat};
				\addlegendentry{$q_y$}
				\addplot[orange,domain=0.02:0.007, dashed]{x*x*x*x*(6000)};
				\addlegendentry{fourth order}
				\addplot[orange,domain=0.02:0.007]{x*x*x*x*x*(300000)};
				\addlegendentry{fifth order}
			\end{axis}
		\end{tikzpicture}}
	\caption{%
		Unsteady vortex from \cref{sec:vortex}: convergence test.
		Left panel: depiction of the initial condition.
		Right panel: error lines, showing that the scheme is indeed of fifth order accuracy.%
	}
	\label{fig:vortex_conv}
\end{figure}

The convergence test is run on Cartesian meshes of sizes
$25^2$, $50^2$, $100^2$, $200^2$, $300^2$, and $400^2$.
The error, denoted by $||\epsilon_h(\bu)||$, is computed as the $\mathbb L^1$ norm of the difference between the approximated solution and the exact one.
\cref{fig:vortex_conv} shows the initial water height for this test case (left panel) and the retrieved fifth order convergence trend expected from theory (right panel).
These results are also reported in \cref{tab:error_table_vortex}, where fifth order accuracy is shown to be achieved.


\begin{table}[!ht]
	\newcolumntype{C}{>{\centering\arraybackslash}p{1.5cm}}
\newcolumntype{R}{>{\raggedleft\arraybackslash}p{1.5cm}}
\newcolumntype{L}{>{\raggedright\arraybackslash}p{1.5cm}}

\pgfplotstableread{vortex_mPDec5_with_orders.dat}\data

\makebox[\textwidth][c]{
	\pgfplotstabletypeset[
		clear infinite,
		empty cells with={---},
		columns={N, h, order-h, hu, order-hu, hv, order-hv},
		columns/N/.style = {int detect, column type = {C}},
		columns/h/.style = {sci zerofill, column type = {C}},
		columns/order-h/.style = {fixed zerofill, column type = {C}},
		columns/hu/.style = {sci zerofill, column type = {C}},
		columns/order-hu/.style = {fixed zerofill, column type = {C}},
		columns/hv/.style = {sci zerofill, column type = {C}},
		columns/order-hv/.style = {fixed zerofill, column type = {C}},
		every head row/.style={
				output empty row,
				before row={
						\toprule%
						& \multicolumn{2}{c}{$h$} & \multicolumn{2}{c}{$q_x$} & \multicolumn{2}{c}{$q_y$} \\
						\cmidrule(lr){2-3}\cmidrule(lr){4-5}\cmidrule(lr){6-7}
						$N_x$ & error & order & error & order & error & order \\
						\cmidrule(lr){1-1}\cmidrule(lr){2-3}\cmidrule(lr){4-5}\cmidrule(lr){6-7}
					}
			},
	]{\data}
}

	\caption{Errors and orders of accuracy, with respect to the number $N_x$ of cells, for the traveling vortex from \cref{sec:vortex}. We indeed observe fifth-order accuracy.}
	\label{tab:error_table_vortex}
\end{table}

\subsection{Lake at rest} 
\label{sec:lake_at_rest}


We now focus on showing the capability of the proposed scheme to exactly preserve the lake at rest steady state, governed by \eqref{eq:lake_at_rest}.
We first tackle the exact capture of the steady state in \cref{sec:latr_exact_capturing},
and we then perform a perturbation analysis in \cref{sec:latr_perturbation_analysis}.

\subsubsection{Exact capturing}
\label{sec:latr_exact_capturing}

In this section, we demonstrate that the proposed scheme is able to exactly capture the lake at rest steady solution in two situations: a fully wet case, where the water height never vanishes, and a wet-dry case, where the water height may vanish.

\subsubsection*{Wet lake at rest}


We first consider the lake at rest steady state given by
\begin{equation*}
	b(x,y) = 0.1\,\sin\left(2\,\pi\,x\right)\,\cos\left(2\,\pi\,y\right),\qquad h(x,y,t) = 1\,-\,b(x,y),\qquad u=v=0,
\end{equation*}
on the computational domain $\Omega\coloneqq [0,1]\times[0,1]$ with periodic boundary conditions and final time $T_f\coloneqq 0.1$.
In this case, we test the scheme with and without the WB modification.
The purpose of this test is twofold.
First, with the WB modification, the lake at rest should be exactly preserved (up to machine precision).
This will confirm that the WB property is satisfied in this case.
Second, without the WB modification, the method should converge to fifth order accuracy.
This will verify the correct implementation of the source term.
Both convergence trends are presented in \cref{fig:larwetconv} and the expected results are obtained.
The resolutions of the Cartesian meshes used for this test are $25^2$, $50^2$, $100^2$ and $200^2$.
%


\pgfplotsset{
  log x ticks with fixed point/.style={
      xticklabel={
        \pgfkeys{/pgf/fpu=true}
        \pgfmathparse{exp(\tick)}%
        \pgfmathprintnumber[fixed relative, precision=3]{\pgfmathresult}
        \pgfkeys{/pgf/fpu=false}
      }
  },
}

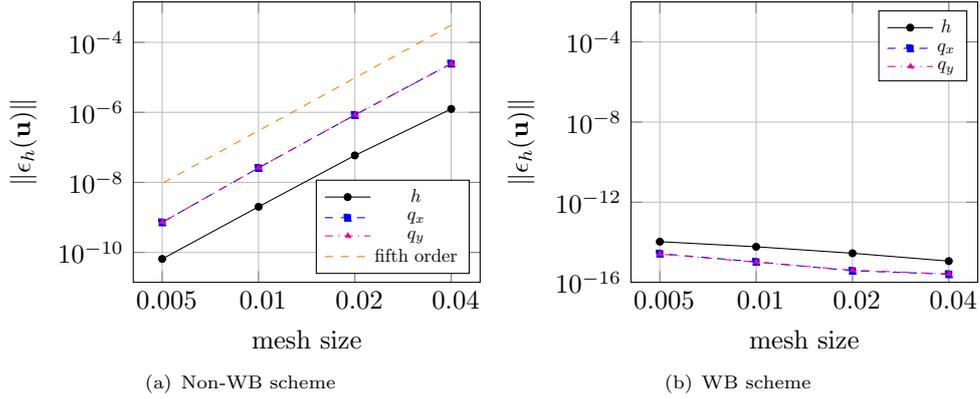
\begin{figure}
	\centering
	\subfigure[Non-WB scheme]{
		\begin{tikzpicture}
			\begin{axis}[
					xmode=log, ymode=log,
                    xtick={0.04,0.02,0.01,0.005},
                    log x ticks with fixed point,
					grid=major,
					xlabel={mesh size},
					ylabel={$\|\epsilon_h(\mathbf{u})\|$},
					xlabel shift = 1 pt,
					ylabel shift = 1 pt,
					legend pos= south east,
					legend style={nodes={scale=0.7, transform shape}},
					width=.47\textwidth
				]
				\addplot[mark=otimes*,mark size=1.3pt,black] table [y=h_NWB, x=size]{LarWet_mPDec5.dat};
				\addlegendentry{$h$}

				\addplot[mark=square*,mark size=1.3pt,blue,dashed] table [y=hu_NWB, x=size]{LarWet_mPDec5.dat};
				\addlegendentry{$q_x$}

				\addplot[mark=triangle*,mark size=1.3pt,magenta,dashdotted] table [y=hv_NWB, x=size]{LarWet_mPDec5.dat};
				\addlegendentry{$q_y$}
				\addplot[orange,domain=0.04:0.005, dashed]{x*x*x*x*x*(3000)};
				\addlegendentry{fifth order}
			\end{axis}
		\end{tikzpicture}}\hfill
	\subfigure[WB scheme]{
		\begin{tikzpicture}
			\begin{axis}[
					xmode=log, ymode=log,
					ymin=1.e-16,ymax=1.e-2,
                    xtick={0.04,0.02,0.01,0.005},
                    log x ticks with fixed point,
					grid=major,
					xlabel={mesh size},
					ylabel={$\|\epsilon_h(\mathbf{u})\|$},
					xlabel shift = 1 pt,
					ylabel shift = 1 pt,
					legend pos= north east,
					legend style={nodes={scale=0.7, transform shape}},
					width=.47\textwidth
				]
				\addplot[mark=otimes*,mark size=1.3pt,black] table [y=h_WB, x=size]{LarWet_mPDec5.dat};
				\addlegendentry{$h$}

				\addplot[mark=square*,mark size=1.3pt,blue,dashed] table [y=hu_WB, x=size]{LarWet_mPDec5.dat};
				\addlegendentry{$q_x$}

				\addplot[mark=triangle*,mark size=1.3pt,magenta,dashdotted] table [y=hv_WB, x=size]{LarWet_mPDec5.dat};
				\addlegendentry{$q_y$}
			\end{axis}
		\end{tikzpicture}}
	\caption{Wet lake at rest from \cref{sec:lake_at_rest}: convergence test. Left panel: non-WB scheme; we observe fifth order accuracy. Right panel: WB scheme; we observe an accuracy up to machine precision for each mesh size.}
    \label{fig:larwetconv}
\end{figure}

\subsubsection*{Wet-dry lake at rest}


We now present a numerical experiment to show the preservation of a lake at rest steady state in the presence of dry areas.
That is to say, the water height will vanish in some parts of the domain.
In particular, by virtue of the mPDeC approach, the proposed method is able to deal with such dry states while having a much relaxed CFL constraint compared to traditional high-order techniques.
Indeed, we can set $\text{CFL} \simeq 1$ rather than $\text{CFL} \simeq 1/12$.

\begin{figure}
	\centering
	\subfigure[Water height $h$]{
		\includegraphics[width=0.50\textwidth,trim={1.2cm 1mm 25mm 1.6cm},clip]{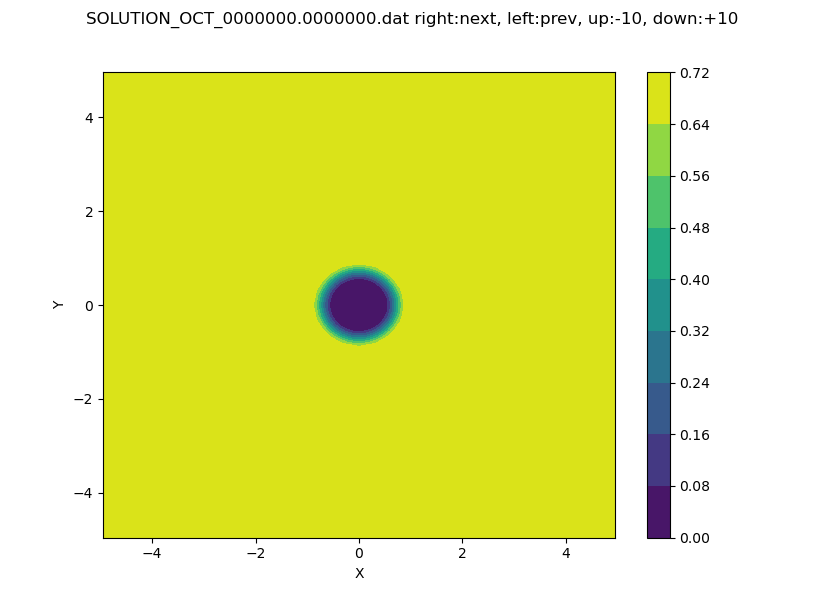}
	}\hfill
	\subfigure[Bathymetry $b$]{
		\includegraphics[width=0.45\textwidth,trim={12cm 1cm 5cm 2.0cm},clip]{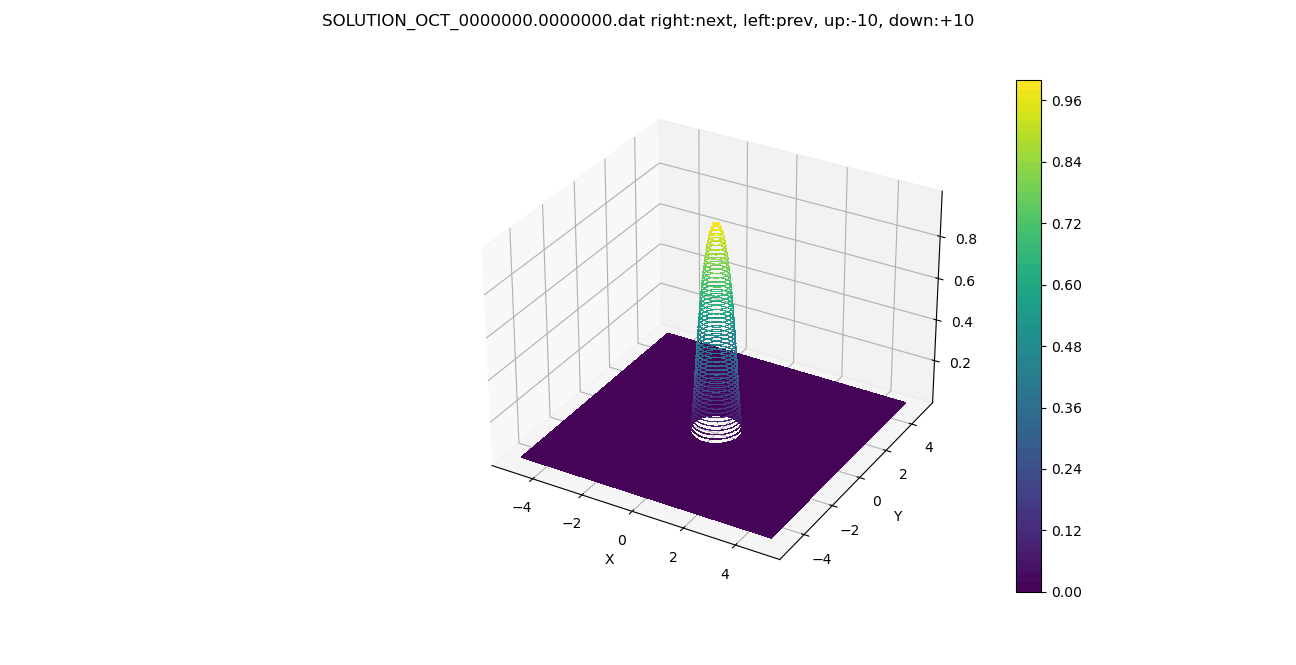}
	}
	\caption{Wet-dry lake at rest from \cref{sec:lake_at_rest}: depiction of the water height (left panel) and of the bathymetry (right panel).}\label{fig:larwetdrysol}
\end{figure}

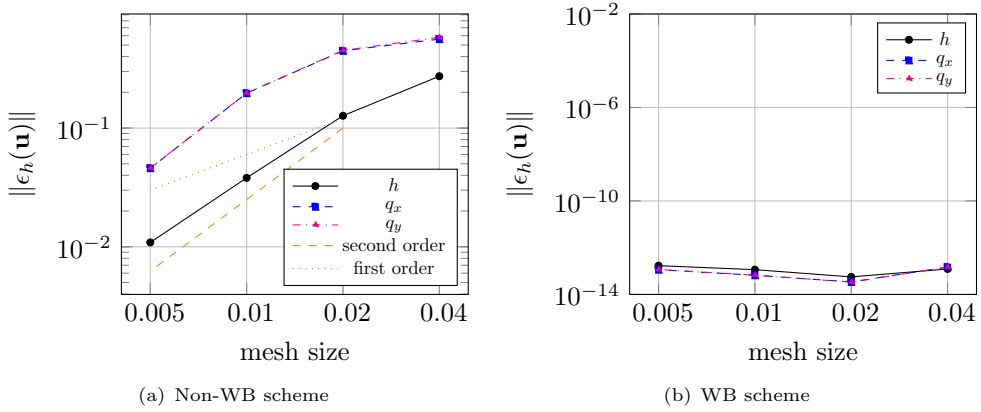
\begin{figure}
	\centering
	\subfigure[Non-WB scheme]{
		\begin{tikzpicture}
			\begin{axis}[
					xmode=log, ymode=log,
                    xtick={0.04,0.02,0.01,0.005},
                    log x ticks with fixed point,
					grid=major,
					xlabel={mesh size},
					ylabel={$\|\epsilon_h(\mathbf{u})\|$},
					xlabel shift = 1 pt,
					ylabel shift = 1 pt,
					legend pos= south east,
					legend style={nodes={scale=0.7, transform shape}},
					width=.47\textwidth
				]
				\addplot[mark=otimes*,mark size=1.3pt,black] table [y=h_NWB, x=size]{LarWetDry_mPDec5.dat};
				\addlegendentry{$h$}

				\addplot[mark=square*,mark size=1.3pt,blue,dashed] table [y=hu_NWB, x=size]{LarWetDry_mPDec5.dat};
				\addlegendentry{$q_x$}

				\addplot[mark=triangle*,mark size=1.3pt,magenta,dashdotted] table [y=hv_NWB, x=size]{LarWetDry_mPDec5.dat};
				\addlegendentry{$q_y$}
				\addplot[orange,domain=0.02:0.005, dashed]{x*x*(250)};
				\addlegendentry{second order}
				\addplot[olive,domain=0.02:0.005, dotted]{x*(6)};
				\addlegendentry{first order}
			\end{axis}
		\end{tikzpicture}}\hfill
	\subfigure[WB scheme]{
		\begin{tikzpicture}
			\begin{axis}[
					xmode=log, ymode=log,
					ymin=1.e-14,ymax=1.e-2,
                    xtick={0.04,0.02,0.01,0.005},
                    log x ticks with fixed point,
					grid=major,
					xlabel={mesh size},
					ylabel={$\|\epsilon_h(\mathbf{u})\|$},
					xlabel shift = 1 pt,
					ylabel shift = 1 pt,
					legend pos= north east,
					legend style={nodes={scale=0.7, transform shape}},
					width=.47\textwidth
				]
				\addplot[mark=otimes*,mark size=1.3pt,black] table [y=h_WB, x=size]{LarWetDry_mPDec5.dat};
				\addlegendentry{$h$}

				\addplot[mark=square*,mark size=1.3pt,blue,dashed] table [y=hu_WB, x=size]{LarWetDry_mPDec5.dat};
				\addlegendentry{$q_x$}

				\addplot[mark=triangle*,mark size=1.3pt,magenta,dashdotted] table [y=hv_WB, x=size]{LarWetDry_mPDec5.dat};
				\addlegendentry{$q_y$}
			\end{axis}
		\end{tikzpicture}}
        \caption{Wet-dry lake at rest from \cref{sec:lake_at_rest}: convergence test. Left panel: non-WB scheme; we observe second order accuracy. Right panel: WB scheme; we observe an accuracy up to machine precision for each mesh size.}
        \label{fig:latrwd}
\end{figure}

We consider, on the domain $\Omega\coloneqq [-5, 5]\times[-5, 5]$, the following bathymetry
\begin{equation}\label{eq:bathymetry_pert}
	b(x,y) \coloneqq
    \begin{cases}
        \exp\left(1-\dfrac{1}{1-r^2}\right), & \text{ if } r^2<1,
        \\ 0,  & \text{ otherwise,} \vphantom{\dfrac 1 2}
    \end{cases}
    \qquad\text{ where } r^2=x^2+y^2.
\end{equation}
This bathymetry represents an island located in the center of the domain.
The water height is defined as $h(x,y,t) \coloneqq  \max(0.7-b(x,y),0)$.
In \cref{fig:larwetdrysol}, we display the water height (left panel) and then bathymetry (right panel). We observe that dry areas occur in the center of the domain, where the island is located.
The test is performed with periodic boundary conditions and final time $T_f\coloneqq 1.$
Just like before, we present the results of a convergence analysis obtained with and without the WB blending.
It should be noted that, due to the discontinuity in the derivative of the water height, the non-WB scheme can achieve
at most second order convergence, while machine precision is expected by the WB version.
The results are reported in \cref{fig:latrwd}; they agree with the expected behavior.

\subsubsection{Perturbation analysis}\label{sec:latr_perturbation_analysis}

Let us consider the computational domain $\Omega\coloneqq [-5,5]\times[-2,2]$, the bathymetry $b$ defined in Equation~\eqref{eq:bathymetry_pert}, and the lake at rest steady state characterized by a total water height $\eta_0\coloneqq 1.5$.
Then, we consider the following perturbation of the steady condition
\begin{equation*}\label{eq:pert_latr}
	\eta \coloneqq  \eta_0+
    \begin{cases}
		0.05 \exp\left(1-\dfrac{1}{(1-\rho^2)^2}\right), & \text{ if } \rho^2<1, \\
		0,                                 & \text{  otherwise,} \vphantom{\dfrac 1 2}
	\end{cases},
\end{equation*}
where we have set
\begin{equation*}
    \rho^2 =9((x+2)^2+(x-0.5)^2).
\end{equation*}
We adopt a Cartesian mesh of $100\times 30$ elements, with periodic boundary conditions, $\text{CFL}=0.8$ and a final time $T_f\coloneqq 0.375$.
We still test the non-WB and WB versions of the method, in order to highlight the advantages of the latter setting. However, in this case, we suppress the blending by adopting all the coefficients $\theta$ of \cref{subsec_Well_balance} equal to 0 in such a way to always use the WB discretization.

\begin{figure}
	\centering
	\subfigure[$t=0$]{\includegraphics[width=0.49\textwidth,trim={50mm 4cm 5cm 45mm},clip]{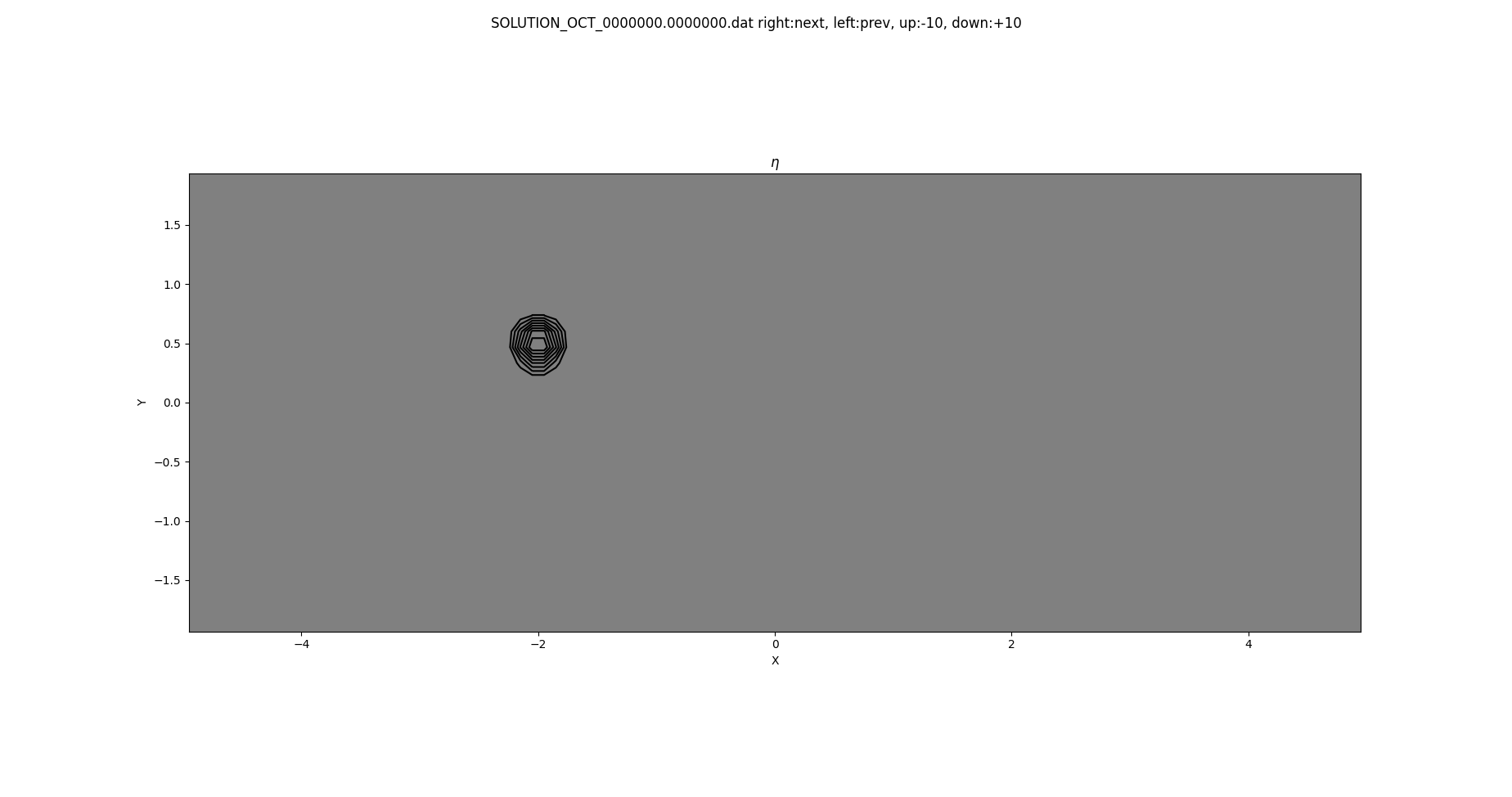}}\\
	\subfigure[$t=0.125$]{\includegraphics[width=0.49\textwidth,trim={50mm 4cm 5cm 45mm},clip]{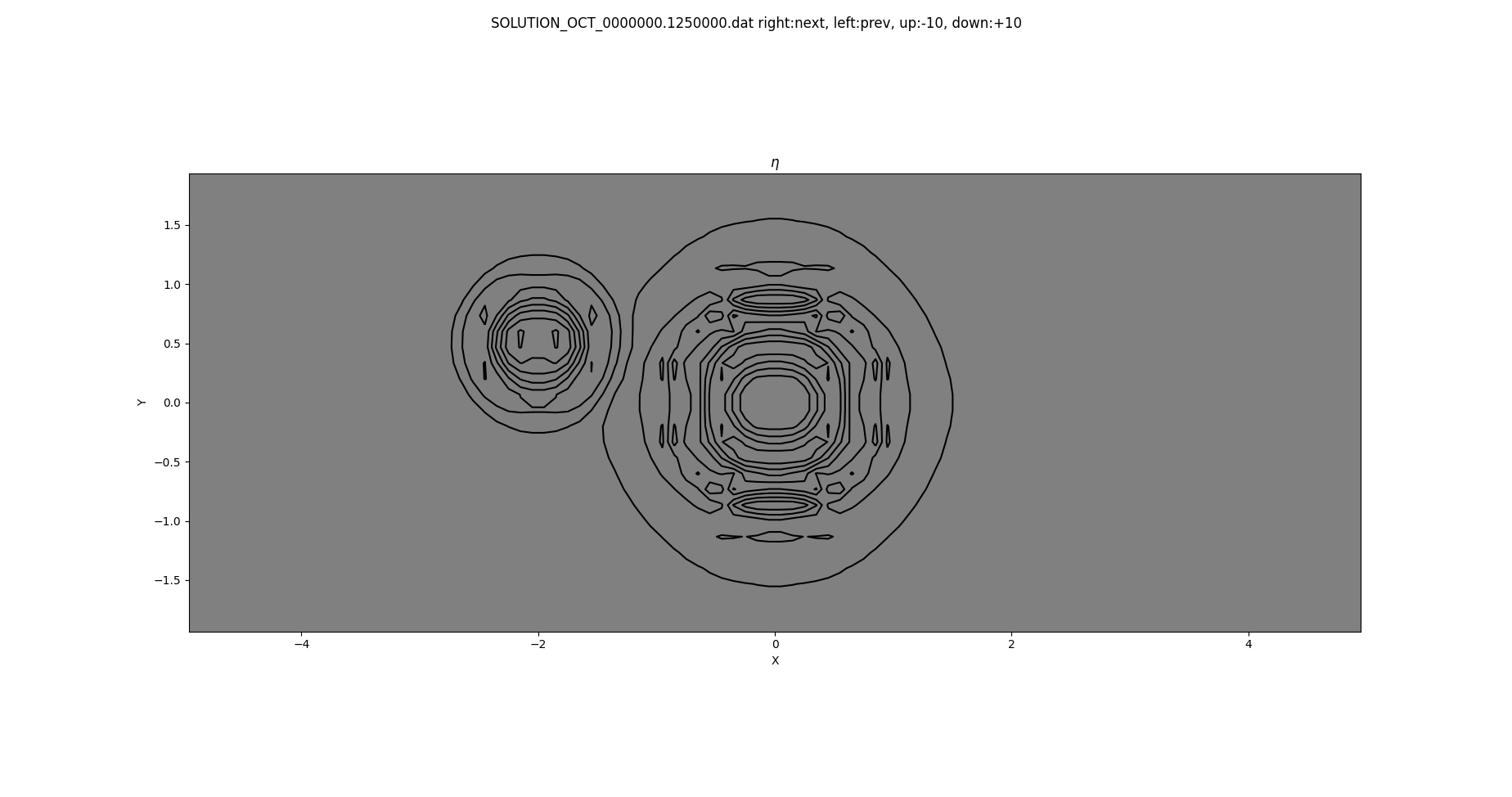}}
	\hfill
	\subfigure[$t=0.125$]{\includegraphics[width=0.49\textwidth,trim={50mm 4cm 5cm 45mm},clip]{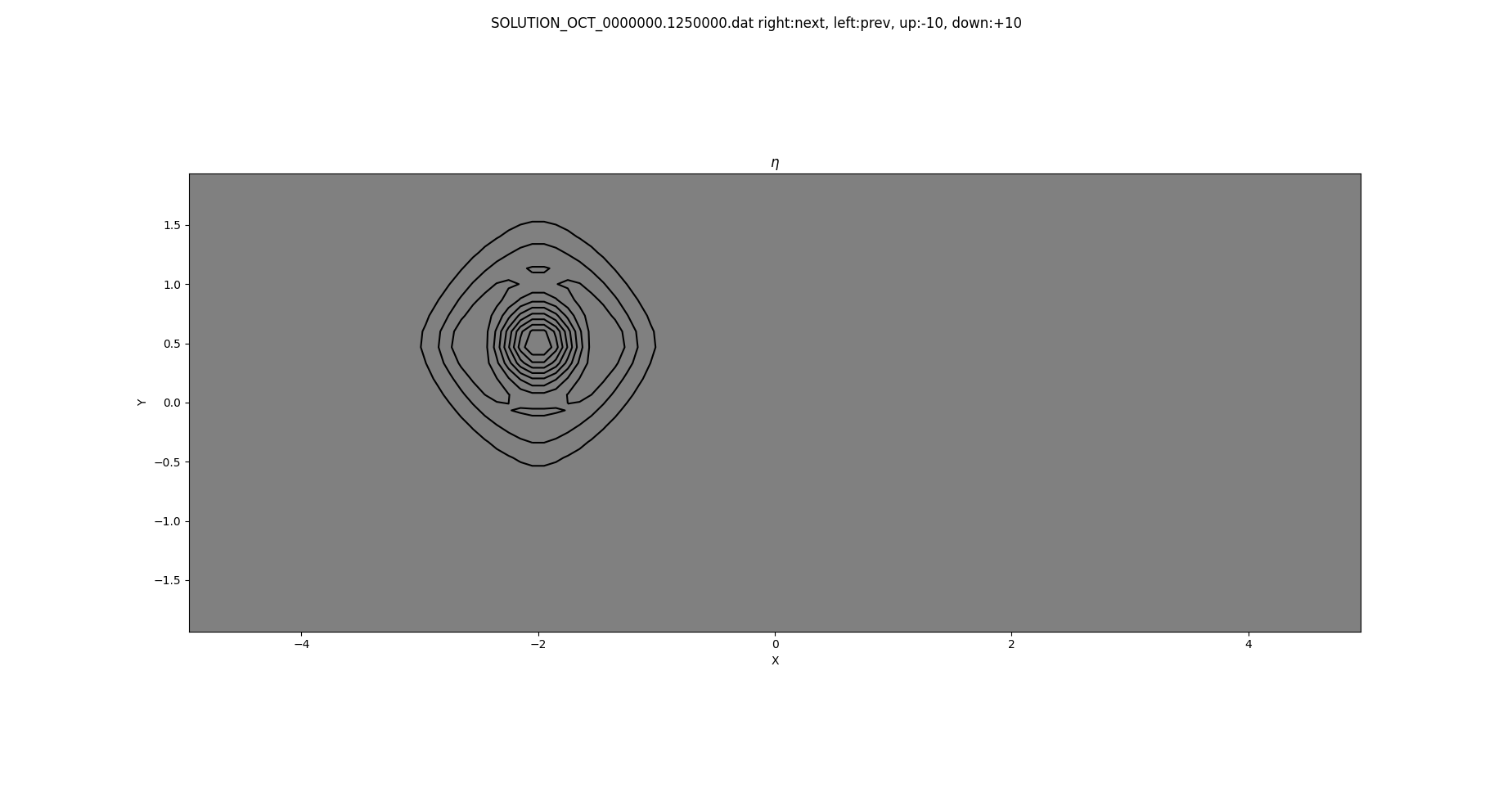}}\\
	\subfigure[$t=0.25$]{\includegraphics[width=0.49\textwidth,trim={50mm 4cm 5cm 45mm},clip]{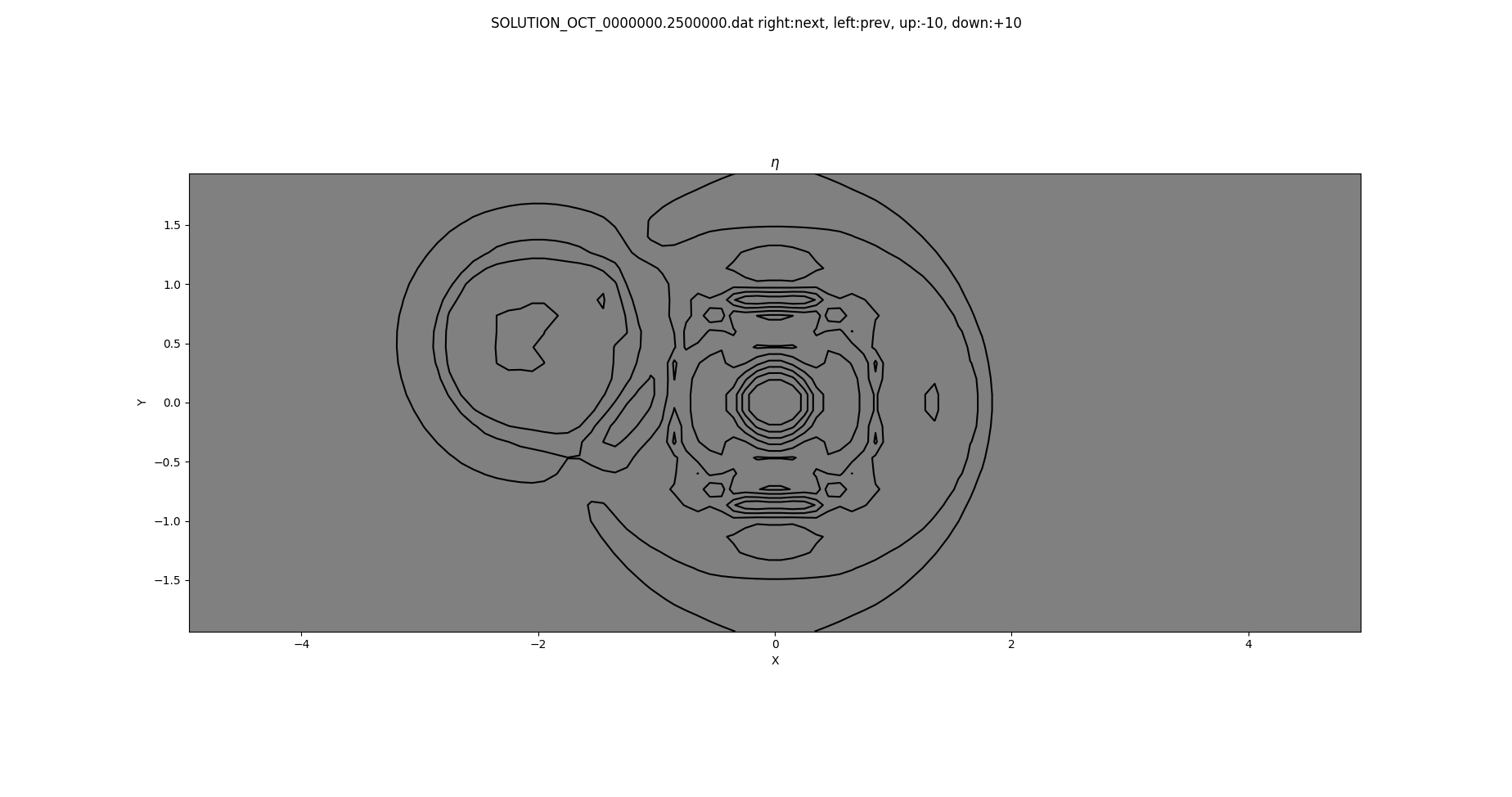}}
	\hfill
	\subfigure[$t=0.25$]{\includegraphics[width=0.49\textwidth,trim={50mm 4cm 5cm 45mm},clip]{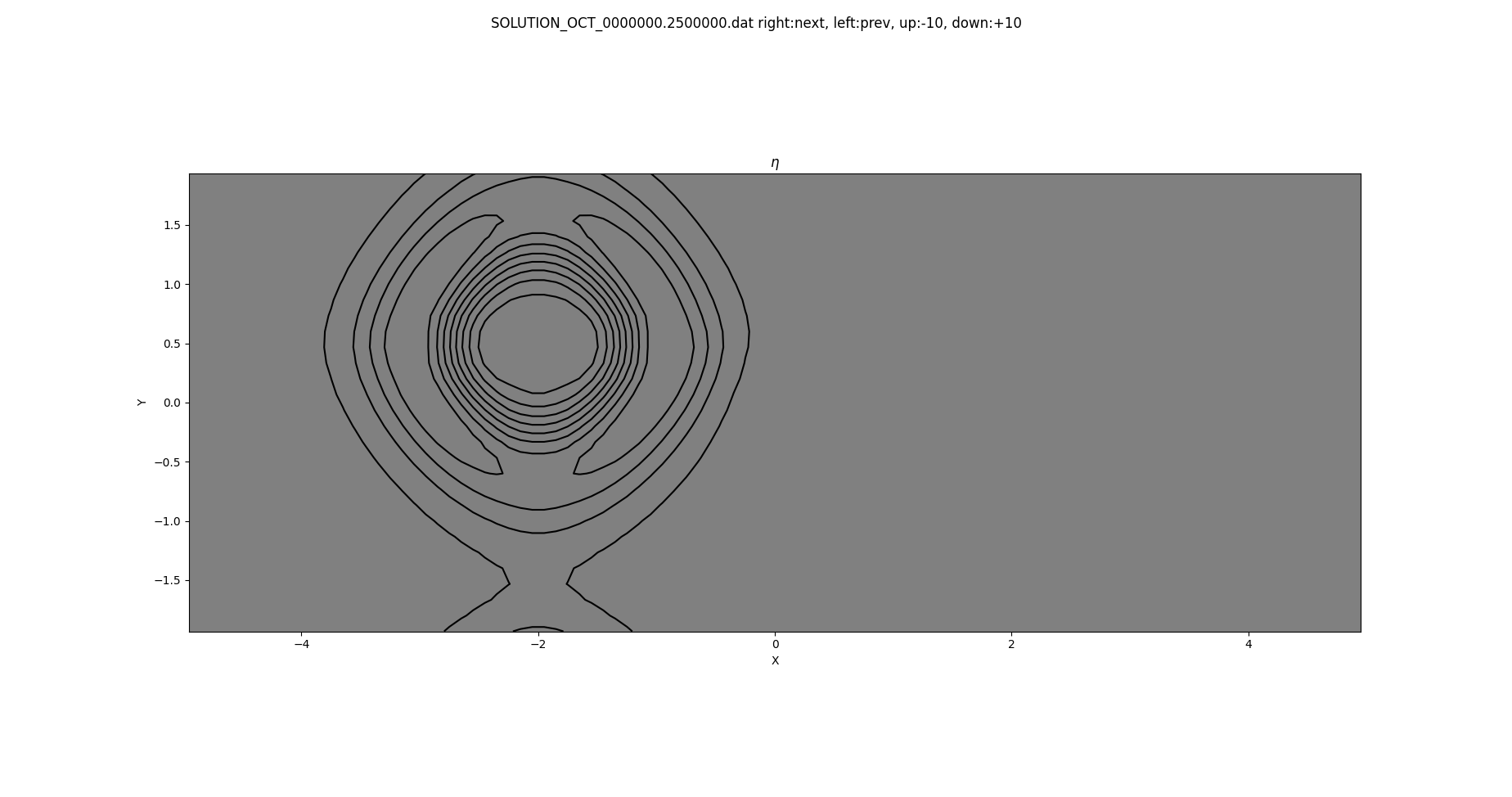}}\\
	\subfigure[$t=0.375$]{\includegraphics[width=0.49\textwidth,trim={50mm 4cm 5cm 45mm},clip]{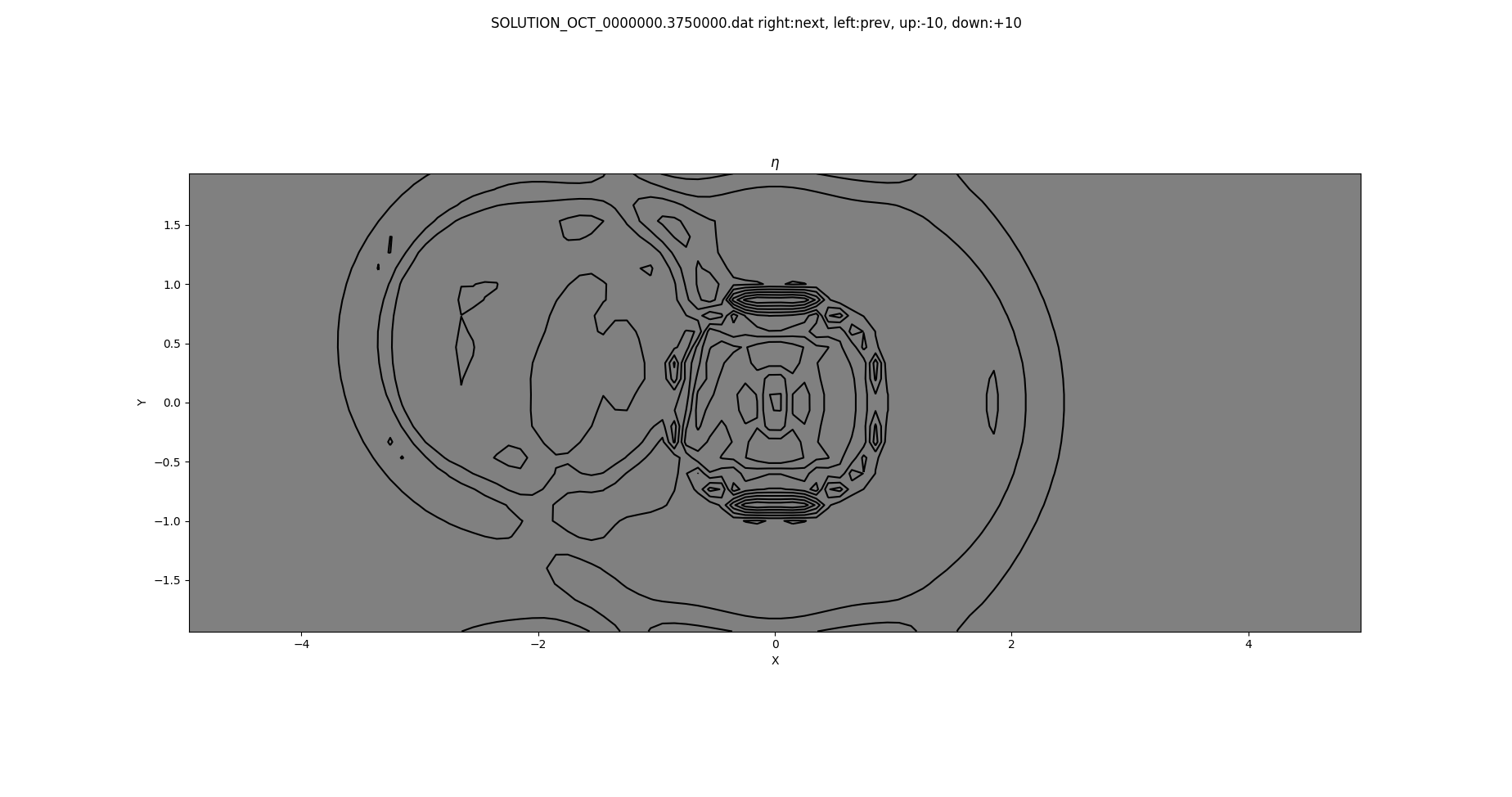}}
	\hfill
	\subfigure[$t=0.375$]{\includegraphics[width=0.49\textwidth,trim={50mm 4cm 5cm 45mm},clip]{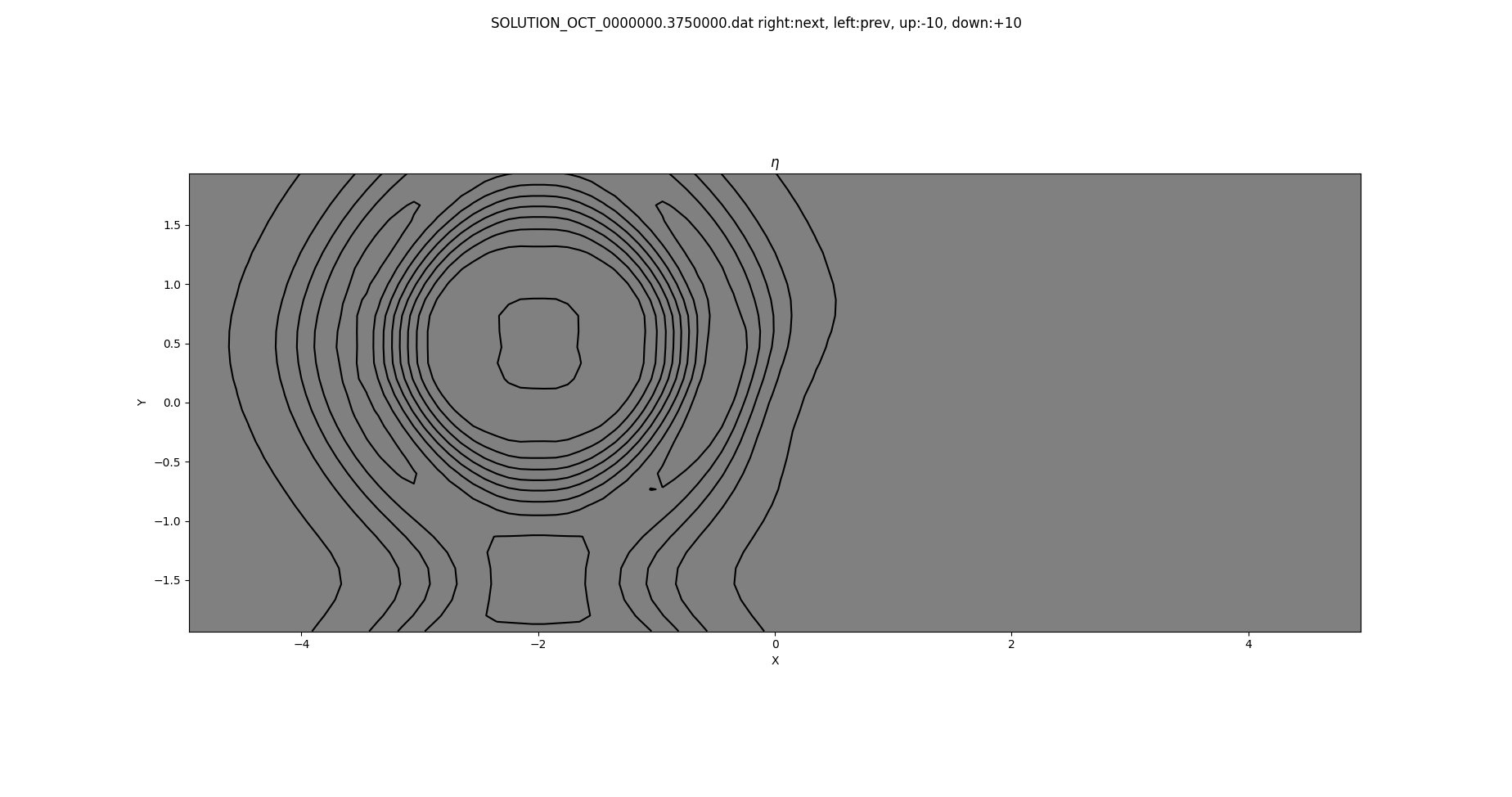}}
	\caption{Perturbation analysis of the lake at rest solution: $\eta=h+b$ isocontours at different times. Top panel, subfigure (a): initial condition. Left panels, subfigures (b), (d) and (f): non-WB scheme; right panels, subfigures (c), (e) and (g): WB scheme.}\label{fig:pert_latr}
\end{figure}

The results are displayed in \cref{fig:pert_latr}.
It can be noticed that the evolution of the perturbation is sharply captured by the WB version of the scheme.
Instead, in the non-WB case, numerical oscillations, due to the discretization error, propagate from the bathymetry
and prevent the proper capturing of the perturbation evolution.

\subsection{Moving equilibria}
\label{sec:simu_moving}


In this section, we test the WB properties of the scheme to capture moving equilibria satisfying \eqref{eq:equilibriumODE}.
As already specified, they are pseudo-1D states.
Therefore, in the context of this section, we focus on the variable $s=x$, and we drop the dependency on $y$, being clear that all quantities are constant along the $y$-direction.
At the numerical level, the variable $y$ does not play any role either.
Hence, the adopted mesh configurations will be characterized by a uniform distribution of cells along the $x$-direction, with various values of $N_x$ ranging from 25 to 200, and a constant number $N_y\coloneqq 5$ of cells along the $y$-direction, with periodic boundary conditions assumed in such direction.

In this 1D frictionless case, moving equilibria are characterized by constant equilibrium variables \eqref{eq:equilibriumvariables}.
Therefore, although there is no closed-form expression of such steady solutions,
they can be computed pointwise, for a given bathymetry,
by solving a cubic equation derived from~\eqref{eq:equilibriumODE}, see for example \cite{delestre2013swashes,ciallella2023arbitrary,micalizzi2023novel}.
The steady flow regime then depends on the prescribed boundary conditions, and is obtained after a transient phase.
We focus here on subcritical and supercritical flows, numerically obtained with the initial and boundary conditions described in \cref{tab:initial_conditions_moving_equilibria},
where the final time~$T_f$ is chosen such that the simulation reaches the steady state (i.e., to make the time residual vanish).
%
We take the following smooth bathymetry
\begin{equation*}
	b(x) \coloneqq  0.05\sin(x-12.5)\exp(1-(x-12.5)^2),
\end{equation*}
on the computational domain $\Omega\coloneqq [0,25]\times[0,1]$.
The gravity constant is set here to $g\coloneqq 9.812$ as in \cite{ciallella2023arbitrary}.

\begin{table}[!ht]
    \centering
    \begin{tabular}{lccccccc}
        \toprule
        Flow regime &
        $T_f$ & $h(x,0)$ & $q(x,0)$ & $h(0,t)$ & $h(25,t)$ & $q(0,t)$ & $q(25,t)$ \\
        \cmidrule(lr){1-8}
        Subcritical &
        200 & $2 - b(x)$ & 0 & --- & 2 & 4.42 & --- \\
        Supercritical &
        50 & $2 - b(x)$ & 0 & 2 & --- & 24 & --- \\
        \bottomrule
        \end{tabular}
    \caption{Initial and boundary conditions for the subcritical and supercritical flows from \cref{sec:simu_moving}. The simulation is always initialized to a lake at rest. Empty cells correspond to Neumann boundary conditions.}
    \label{tab:initial_conditions_moving_equilibria}
\end{table}

Again, we test the WB and the non-WB versions of the scheme.
We emphasize that we do not try to merely exactly preserve the steady solution, but to capture it: the WB numerical scheme is expected to converge towards the steady solution with machine accuracy, even after the transient, unsteady phase.

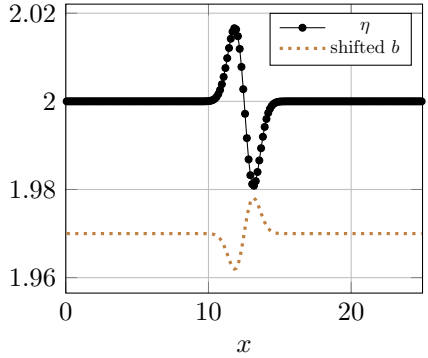
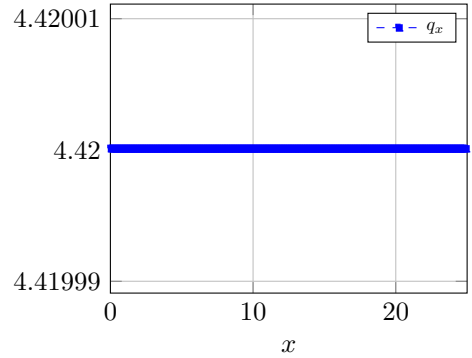
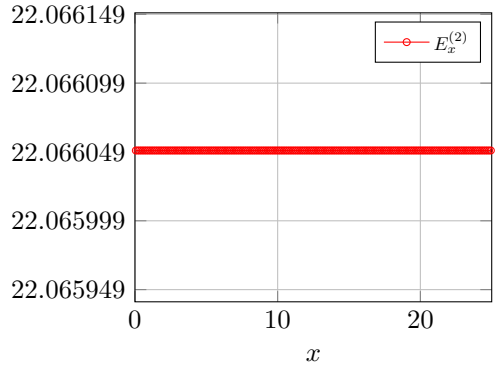
\begin{figure}
	\centering
	\subfigure[free surface $\eta$ and bathymetry $b$, shifted and rescaled]{
		\begin{tikzpicture}
			\begin{axis}[
					xmin=0,xmax=25,
					grid=major,
					xlabel={$x$},
					xlabel shift = 1 pt,
					ylabel shift = 1 pt,
					legend pos= north east,
					legend style={nodes={scale=0.7, transform shape}},
					width=.48\textwidth
				]
				\addplot[mark=otimes*,mark size=1.3pt,black] table [y expr=(\thisrow{h_WB}+\thisrow{bath}), x=size]{Subcrit_mPDec5_solution.dat};
				\addlegendentry{$\eta$}
				\addplot[dotted,very thick,brown] table [y expr=(\thisrow{bath})*0.15+1.97, x=size]{Subcrit_mPDec5_solution.dat};
				\addlegendentry{shifted $b$}
			\end{axis}
		\end{tikzpicture}}\hfill
	\subfigure[discharge $q_x$]{
		\begin{tikzpicture}
			\begin{axis}[
					xmin=0,xmax=25,
					grid=major,
					ymin=4.419989, ymax=4.420011,
					xlabel={$x$},
					yticklabel style={/pgf/number format/.cd,fixed,precision=6},
					xlabel shift = 1 pt,
					ylabel shift = 1 pt,
					legend pos= north east,
					legend style={nodes={scale=0.7, transform shape}},
					width=.48\textwidth
				]
				\addplot[mark=square*,mark size=1.3pt,blue,dashed] table [y=hu_WB, x=size]{Subcrit_mPDec5_solution.dat};
				\addlegendentry{$q_x$}
			\end{axis}
		\end{tikzpicture}}
    \subfigure[second equilibrium variable $E_x^{(2)}$]{
		\begin{tikzpicture}
			\begin{axis}[
					xmin=0,xmax=25,
					grid=major,
					ymin=22.06594, ymax=22.06615,
					xlabel={$x$},
					yticklabel style={/pgf/number format/.cd,fixed,precision=6},
					xlabel shift = 1 pt,
					ylabel shift = 1 pt,
					legend pos= north east,
					legend style={nodes={scale=0.7, transform shape}},
					width=.48\textwidth
				]
				\addplot[mark=o,mark size=1.3pt,red] table [y=E, x=size]{Subcrit_mPDec5_solution_out.dat};
				\addlegendentry{$E_x^{(2)}$}
			\end{axis}
		\end{tikzpicture}}
    \caption{Subcritical flow, test case from \cref{sec:simu_moving}.
    Top left panel: free surface water level $\eta$, and bathymetry $b$ rescaled by a factor of $0.15$ and shifted by $1.97$.
    Top right panel: discharge $q_x$.
    Bottom panel: second equilibrium variable $E_x^{(2)}$.}\label{fig:subsol}
\end{figure}
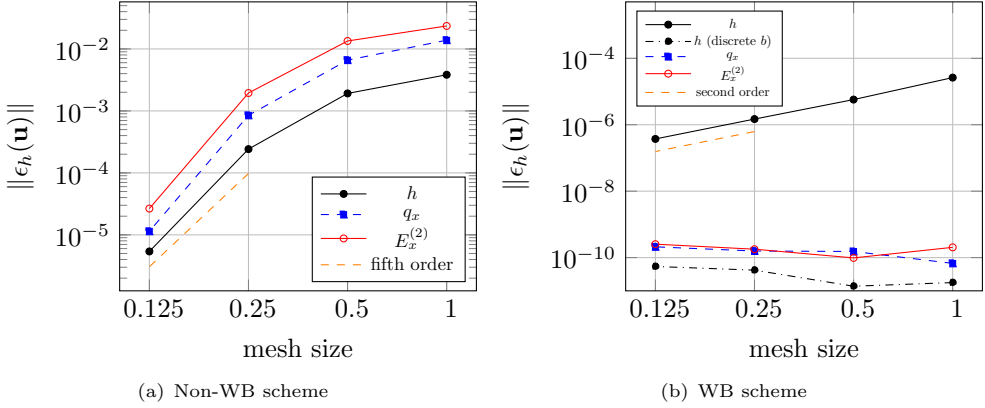
\begin{figure}
	\centering
	\subfigure[Non-WB scheme]{
		\begin{tikzpicture}
			\begin{axis}[
					xmode=log, ymode=log,
                    xtick={1,0.5,0.25,0.125},
                    log x ticks with fixed point,
					grid=major,
					xlabel={mesh size},
					ylabel={$\|\epsilon_h(\mathbf{u})\|$},
					xlabel shift = 1 pt,
					ylabel shift = 1 pt,
					legend pos= south east,
					legend style={nodes={scale=0.7, transform shape}},
					width=.48\textwidth
				]
				\addplot[mark=otimes*,mark size=1.3pt,black] table [y=h_NWB, x=size]{Subcrit_mPDec5.dat};
				\addlegendentry{$h$}

				\addplot[mark=square*,mark size=1.3pt,blue,dashed] table [y=hu_NWB, x=size]{Subcrit_mPDec5.dat};
				\addlegendentry{$q_x$}
				\addplot[mark=o,mark size=1.3pt,red] table [y=E_NWB, x=size]{Subcrit_mPDec5.dat};
				\addlegendentry{$E_x^{(2)}$}

				\addplot[orange,domain=0.25:0.125, dashed]{x*x*x*x*x*(0.1)};
				\addlegendentry{fifth order}
			\end{axis}
		\end{tikzpicture}}\hfill
	\subfigure[WB scheme]{
		\begin{tikzpicture}
			\begin{axis}[
					xmode=log, ymode=log,
					ymin=1.e-11, ymax=5.e-3,
                    xtick={1,0.5,0.25,0.125},
                    log x ticks with fixed point,
					grid=major,
					xlabel={mesh size},
					ylabel={$\|\epsilon_h(\mathbf{u})\|$},
					xlabel shift = 1 pt,
					ylabel shift = 1 pt,
					legend pos= north west,
					legend style={nodes={scale=0.5, transform shape}},
					width=.48\textwidth
				]
				\addplot[mark=otimes*,mark size=1.3pt,black] table [y=h_WB, x=size]{Subcrit_mPDec5.dat};
				\addlegendentry{$h$}

				\addplot[mark=otimes*,mark size=1.3pt,black,dashdotted] table [y=h_discreteWB, x=size]{Subcrit_mPDec5.dat};
				\addlegendentry{$h$ (discrete $b$)}

				\addplot[mark=square*,mark size=1.3pt,blue,dashed] table [y=hu_WB, x=size]{Subcrit_mPDec5.dat};
				\addlegendentry{$q_x$}
				\addplot[mark=o,mark size=1.3pt,red] table [y=E_WB, x=size]{Subcrit_mPDec5.dat};
				\addlegendentry{$E_x^{(2)}$}

				\addplot[orange,domain=0.25:0.125, dashed]{x*x*(0.00001)};
				\addlegendentry{second order}
			\end{axis}
		\end{tikzpicture}}
    \caption{Subcritical flow, test case from \cref{sec:simu_moving}: convergence test. Left panel: non-WB scheme; we observe fifth-order accuracy. Right panel: WB scheme; we observe machine precision accuracy for the equilibrium variables $q_x$ and $E_x^{(2)}$ and $h$ when considering the discrete bathymetry, and second-order accuracy for $h$ when considering the exact bathymetry function.}\label{fig:subconv}
\end{figure}

We start by presenting the numerical results obtained for the subcritical steady flow.
The solution computed with the WB method, with $N_x\coloneqq 200$, is presented in \cref{fig:subsol}.
We display three quantities: the water height $h$, the $x$-discharge $q_x$, and the second component of the equilibrium variables $E_x^{(2)}\coloneqq \frac 1 2 \frac{q_x^2}{h^2} + g(h+b)$.
Recall that both $q_x$ and $E_x^{(2)}$ should be constant in this case;
we can indeed appreciate the ability of the WB blending to capture constant $q_x$ and $E_x^{(2)}$.
The exact capture of $q_x$ and~$E_x^{(2)}$ is also visible from the results of the convergence test reported in \cref{fig:subconv}.
In particular, we observe that the WB version of the scheme is able to obtain machine precision errors with respect to such variables.
%
Notice that the errors on the water height can be computed following two approaches: the first one, which takes as a reference the exact bathymetry function $b(x)$; the second one, which considers the discrete bathymetry in cell average $b_{i,j}$. 
Usually the first approach is employed for classical convergence analysis, however the second one is very common in the field of 
well-balanced schemes to check whether the scheme is able to preserve the discrete version of the considered equilibrium.
While the first method computes the error by using the exact bathymetry evaluated at quadrature points, the second one considers the 
reconstructed bathymetry to measure the reference equilibrium.
For the first convergence test, in line with \cref{rmk:WBLO}, we expect the error $h$ to scale with first order.
However, a second order superconvergence is obtained, due to the exact preservation of~$q_x$ and $E_x^{(2)}$.
For the second convergence analysis, the error $h$ (discrete $b$) provides the proof of the exact preservation of 
discrete steady states of the WB scheme with machine precision obtained for all meshes.
On the other hand, the non-WB scheme produces, as expected, bigger errors which scale with the expected fifth order.
Let us notice that a very high level of mesh refinement would be needed in order to obtain, with the non-WB scheme, errors comparable to the ones obtained with the WB version, especially on $q_x$ and $E_x^{(2)}$.
This ensures that, for a given error, the WB method has a much smaller computational cost than the non-WB one.

Similar considerations apply to the supercritical case. The numerical solution computed with the WB method, with $N_x\coloneqq 200$, is reported in \cref{fig:supersol}. In addition, the convergence plots of both WB and non-WB schemes can be found in \cref{fig:superconv}.
Also in this case, we observe the same features and trends as before:
the ability of the WB version to capture $h$ with discrete bathymetry and the constant equilibrium variables $q_x$ and $E_x^{(2)}$ up to machine precision, and to obtain much smaller errors with respect to the non-WB version.
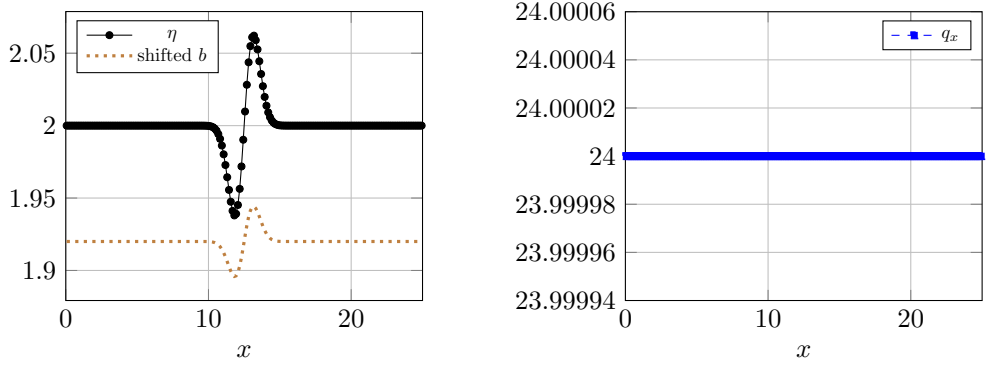
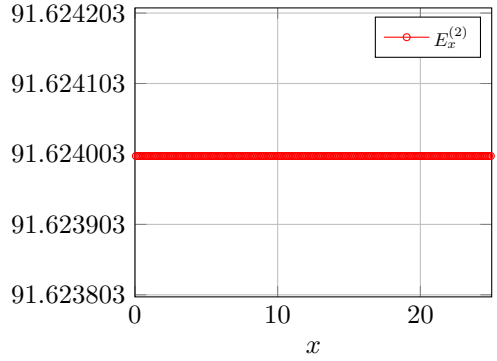
\begin{figure}
	\centering
	\subfigure[free surface $\eta$ and bathymetry $b$, shifted and rescaled]{
		\begin{tikzpicture}
			\begin{axis}[
					grid=major,
					xmin=0,xmax=25,
					legend pos= north west,
					xlabel={$x$},
					yticklabel style={/pgf/number format/.cd,fixed,precision=3},
					legend style={nodes={scale=0.7, transform shape}},
					width=.48\textwidth
				]
				\addplot[mark=otimes*,mark size=1.3pt,black] table [y expr=(\thisrow{Depth}+\thisrow{Bathymetry}), x=CoordinateX]{Supercritical_mPDec5_200-5.dat};
				\addlegendentry{$\eta$}
				\addplot[dotted,very thick,brown] table [y expr=(\thisrow{Bathymetry})*0.45+1.92, x=CoordinateX]{Supercritical_mPDec5_200-5.dat};
				\addlegendentry{shifted $b$}

			\end{axis}
		\end{tikzpicture}}\hfill
	\subfigure[discharge $q_x$]{
		\begin{tikzpicture}
			\begin{axis}[
					grid=major,
					xmin=0,xmax=25,
					ymin=23.99994,
					ymax=24.00006,
					xlabel={$x$},
					yticklabel style={/pgf/number format/.cd,fixed,precision=6},
					ytick={23.99994,23.99996,23.99998,24,24.00002,24.00004,24.00006},
					legend pos= north east,
					legend style={nodes={scale=0.7, transform shape}},
					width=.48\textwidth
				]
				\addplot[mark=square*,mark size=1.3pt,blue,dashed] table [y=qx, x=CoordinateX]{Supercritical_mPDec5_200-5.dat};
				\addlegendentry{$q_x$}

			\end{axis}
		\end{tikzpicture}}
	\subfigure[second equilibrium variable $E_x^{(2)}$]{
		\begin{tikzpicture}
			\begin{axis}[
					grid=major,
					xmin=0,xmax=25,
					ymin=91.62380,
					ymax=91.62421,
					xlabel={$x$},
					yticklabel style={/pgf/number format/.cd,fixed,precision=6},
					legend pos= north east,
					legend style={nodes={scale=0.7, transform shape}},
					width=.48\textwidth
				]
				\addplot[mark=o,mark size=1.3pt,red] table [y=G, x=CoordinateX]{Supercritical_mPDec5_200-5.dat};
				\addlegendentry{$E_x^{(2)}$}

			\end{axis}
		\end{tikzpicture}}
        \caption{Supercritical flow, test case from \cref{sec:simu_moving}.
        Top left panel: free surface water level $\eta$, and bathymetry $b$ rescaled by a factor of $0.45$ and shifted by $1.92$.
        Top right panel: discharge $q_x$.
        Bottom panel: second equilibrium variable $E_x^{(2)}$.}\label{fig:supersol}
\end{figure}
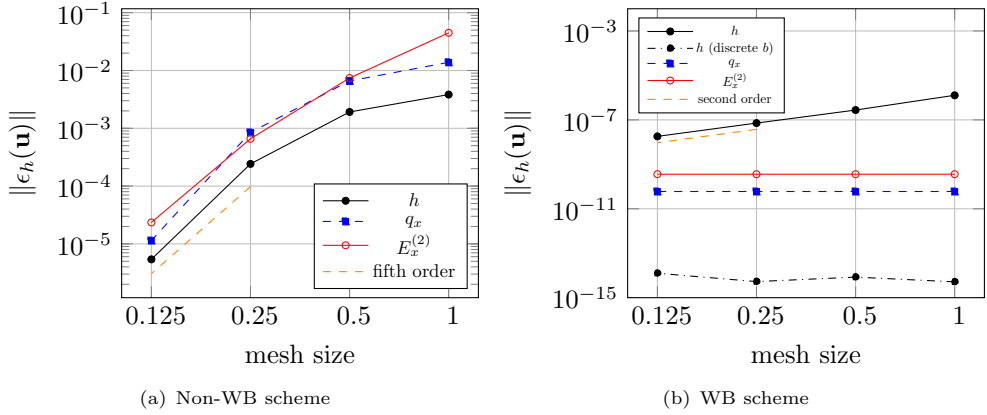
\begin{figure}
	\centering
	\subfigure[Non-WB scheme]{
		\begin{tikzpicture}
			\begin{axis}[
					xmode=log, ymode=log,
                    xtick={1,0.5,0.25,0.125},
                    log x ticks with fixed point,
					grid=major,
					xlabel={mesh size},
					ylabel={$\|\epsilon_h(\mathbf{u})\|$},
					xlabel shift = 1 pt,
					ylabel shift = 1 pt,
					legend pos= south east,
					legend style={nodes={scale=0.7, transform shape}},
					width=.48\textwidth
				]
				\addplot[mark=otimes*,mark size=1.3pt,black] table [y=h_NWB, x=size]{Supercrit_mPDec5.dat};
				\addlegendentry{$h$}

				\addplot[mark=square*,mark size=1.3pt,blue,dashed] table [y=hu_NWB, x=size]{Supercrit_mPDec5.dat};
				\addlegendentry{$q_x$}

				\addplot[mark=o,mark size=1.3pt,red] table [y=E_NWB, x=size]{Supercrit_mPDec5.dat};
				\addlegendentry{$E_x^{(2)}$}

				\addplot[orange,domain=0.25:0.125, dashed]{x*x*x*x*x*(0.1)};
				\addlegendentry{fifth order}
			\end{axis}
		\end{tikzpicture}}\hfill
	\subfigure[WB scheme]{
		\begin{tikzpicture}
			\begin{axis}[
					xmode=log, ymode=log,
					ymin=1.e-15, 
					ymax=1.e-2,
                    			xtick={1,0.5,0.25,0.125},
                   		        log x ticks with fixed point,
					grid=major,
					xlabel={mesh size},
					ylabel={$\|\epsilon_h(\mathbf{u})\|$},
					xlabel shift = 1 pt,
					ylabel shift = 1 pt,
					legend pos= north west,
					legend style={nodes={scale=0.5, transform shape}},
					width=.48\textwidth
				]
				\addplot[mark=otimes*,mark size=1.3pt,black] table [y=h_WB, x=size]{Supercrit_mPDec5.dat};
				\addlegendentry{$h$}

				\addplot[mark=otimes*,mark size=1.3pt,black,dashdotted] table [y=h_discreteWB, x=size]{Supercrit_mPDec5.dat};
				\addlegendentry{$h$ (discrete $b$)}

				\addplot[mark=square*,mark size=1.3pt,blue,dashed] table [y=hu_WB, x=size]{Supercrit_mPDec5.dat};
				\addlegendentry{$q_x$}

				\addplot[mark=o,mark size=1.3pt,red] table [y=E_WB, x=size]{Supercrit_mPDec5.dat};
				\addlegendentry{$E_x^{(2)}$}

				\addplot[orange,domain=0.25:0.125, dashed]{x*x*(0.0000006)};
				\addlegendentry{second order}
			\end{axis}
		\end{tikzpicture}}
    \caption{Supercritical flow, test case from \cref{sec:simu_moving}: convergence test. Left panel: non-WB scheme; we observe fifth-order accuracy. Right panel: WB scheme; we observe machine precision accuracy for the equilibrium variables $q_x$ and $E_x^{(2)}$ and $h$ when considering the discrete bathymetry, and second-order accuracy for $h$ when considering the exact bathymetry function.}\label{fig:superconv}
\end{figure}


\subsection{Flooding simulations}
\label{sec:simu_flooding}

We finally present the numerical results of flooding simulations performed with the proposed high-order WB positivity-preserving method.
While we so far have focused on the validation of the proposed method on standard academic test cases,
we now deal with more challenging applications.
These applications correspond to waves over dry areas, and prove the suitability of the proposed approach in the context of real-world situations.
We start by presenting a wave over a dry island in \cref{sec:wave_over_dry_island}, and then we move to the simulation of a tsunami over three obstacles in \cref{sec:tsunami_on_three_obstacles}.

\subsubsection{Wave over a dry island}
\label{sec:wave_over_dry_island}

In this test, we simulate a wave over a dry island. 
The computational domain is the rectangular region $\Omega\coloneqq [-5, 5]\times[-2, 2]$, partitioned into a mesh with $400 \times 120$ elements.
We refer to \cite{ciallella2022arbitrary}, Section 6.8, for the bathymetry function $b(x,y)$ and the specific initial and boundary conditions.
The simulation was run until a final time $T_f\coloneqq 5$, with a CFL number set to 0.9.

\begin{figure}
	\centering
	\subfigure[$t=0$]{\includegraphics[width=0.48\textwidth, trim={150 80 120 100}, clip]{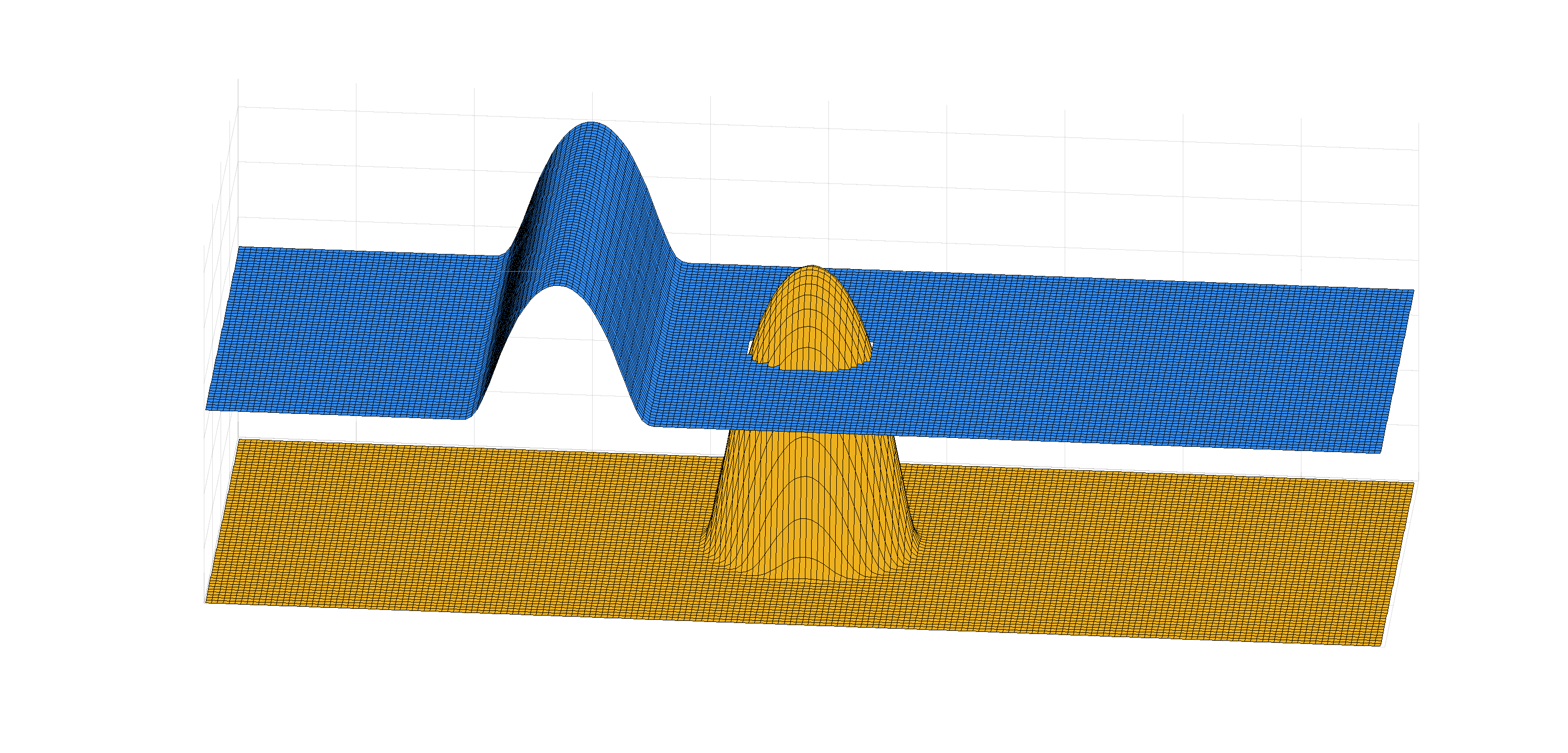}}\\
	\foreach \n in {21, 41, 81, 121, 161, 251} {
			\pgfmathsetmacro{\myresult}{(\n-1)*5.0/250}
			\subfigure[\(t=\pgfmathprintnumber{\myresult}\)]{\includegraphics[width=0.48\textwidth, trim={150 80 120 150}, clip]{island/SOL\n.png}}
		}
	\caption{Wave over a dry island test case from \cref{sec:wave_over_dry_island}: $\eta\coloneqq h+b$ and $b$ at different times.}
	\label{fig:island_simulations}
\end{figure}

The results are presented at various times in \cref{fig:island_simulations}.
The variable $\eta$, along with the bathymetry $b$, have been displayed.
Indeed, it allows for a clearer understanding of the underlying physics.

The simulation starts with a background state moving from left to right at speed $u = 1$, propelling the wave towards the island. This causes the island to get wet from the left side and to dry from the right side.
Thus, the top of the island, initially dry, undergoes multiple wet and dry cycles throughout the whole simulation, without encountering any issue related to negative water height.
This is not guaranteed for classical time integration schemes, among which SSPRK schemes, for such high CFL numbers.
Various structures are observable in this simulation like vortices and shocks, and the recurring wetting/drying processes are optimally tackled by the proposed scheme.

\subsubsection{Tsunami on three obstacles}
\label{sec:tsunami_on_three_obstacles}


Finally, the simulation of a tsunami over several obstacles is presented.
Simulations of this kind are often performed \cite{guermond2022well} since they represent
a good starting point to move towards the simulation of real coastal engineering problems.
In this simulation, we consider a shock impacting three conical obstacles.
More specifically, we consider the domain $\Omega\coloneqq [-5, 7]\times[-2, 2]$, partitioned into $960\times 320$ elements, and the bathymetry
\begin{equation*}
	b(x,y)\coloneqq \sum_{i=1}^3 b_i(x,y)+
    \begin{cases}
		1+0.2 x, \quad    & \text{for } x<0,            \\
		1, \quad          & \text{for } 0\leq x \leq 3, \\
		1+0.4(x-3), \quad & \text{for } x>3,            \\
	\end{cases}
\end{equation*}
with $b_i(x,y)\coloneqq c(x,y,x_i,y_i,R_i,A_i)$, where $c$ is a cone function defined as
\begin{equation*}
	\begin{split}
	&c(x,y,x_c,y_c,R,A)\coloneqq \\
	&\begin{cases}
		\frac A R \left(R-\sqrt{(x-x_c)^2+(y-y_c)^2}\right), \quad & \text{if }\sqrt{(x-x_c)^2+(y-y_c)^2}<R, \\
		0, \quad                                       & \text{otherwise.}                       \\
	\end{cases}
	\end{split}
\end{equation*}
In particular, we have $R_i\coloneqq 0.5$ and $A_i\coloneqq 3$, for all $i$, and $(x_1,y_1)\coloneqq (1,-1)^T$, $(x_2,y_2)\coloneqq (1,1)$ and $(x_3,y_3)\coloneqq (2,0)$.
The initial condition is given by
\begin{equation}
    \label{eq:tsunami_initial_condition}
	\begin{bmatrix}
		h \\
		u \\
		v
	\end{bmatrix}(x,y,0)\coloneqq \begin{cases}
		[1.5-b(x,y),4,0]^T,\quad & \text{if}~x<-3.5, \\
		[0,0,0]^T,\quad          & \text{otherwise}.
	\end{cases}
\end{equation}
The prescribed boundary conditions are
\begin{itemize}[nosep]
    \item inflow at the left of the domain, obtained by imposing $q(x=-5,y,t)\coloneqq 3(1+\cos(2\pi t))e^{-2t}$;
    \item transmissive at the right of the domain;
    \item solid walls at the top and bottom of the domain.
\end{itemize}
We remark that, to simulate a realistic configuration, a time-dependent inlet condition has been chosen to represent a series of waves impacting the obstacles after the tsunami.
The final time is $T_f\coloneqq 3$, and we take a CFL condition of $0.8$ for added stability.

\begin{figure}
	\centering
	\subfigure[$t=0$\label{fig:tsunami_0}]{\includegraphics[width=0.48\textwidth, trim={250 100 210 200}, clip]{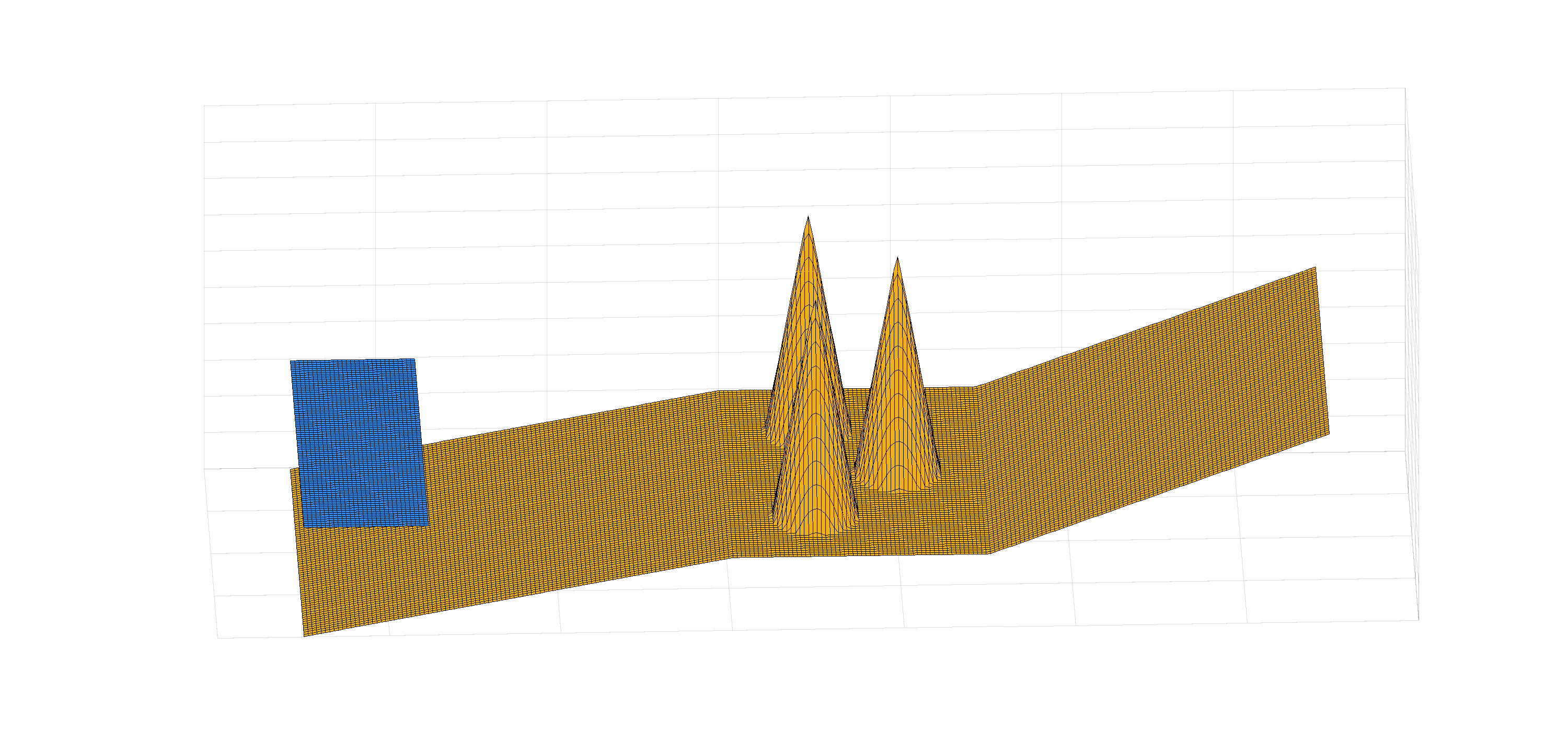}}\\
	\foreach \n in {51, 101, 151, 201, 251, 301} {
			\pgfmathsetmacro{\myresult}{3.0*(\n-1)/300}
			\subfigure[\(t=\myresult\)\label{fig:tsunami_\n}]{\includegraphics[width=0.48\textwidth, trim={250 100 210 200}, clip]{tsunami/SOL\n.png}}
		}
	\caption{Tsunami on three obstacles, test case from \cref{sec:tsunami_on_three_obstacles}: $\eta:h+b$ and $b$ at different times.}
	\label{fig:overall_tsunami}
\end{figure}

The results are reported in \cref{fig:overall_tsunami}.
We start from an initial configuration where the majority of the domain is dry and where the initial tsunami is represented by a discontinuity in the water height, defined in \eqref{eq:tsunami_initial_condition}.
Already from the first snapshots in \cref{fig:tsunami_101,fig:tsunami_151}, we can appreciate the wetting process happening with several structures forming on the right of the three bodies.
Thanks to the time-dependent inlet condition, the dynamic of the simulation keeps evolving with shock interactions occurring due to the crushing between new and old wave fronts,
while wetting and drying processes keep happening in many parts of the domain.

This simulation best represents the potential of this framework, which is able to retain high-order accuracy, preserve important structures of the model, and deal with complex fluid phenomena.
Simulations of this kind are not only challenging but also computationally expensive, due to the accuracy required to capture all flow features.
The choice of the considered time-stepping scheme, able to relax the typical severe CFL constraints imposed by positivity preservation, has a huge impact on the computational resources needed to perform these simulations.
In fact, keeping the same fifth-order accuracy, we are able to consistently reduce the computational time with respect to classical time integration techniques, provably guaranteeing positivity of the discrete water height.
Moreover, the well-balanced procedure is non-intrusive and computationally cheap, and it is able to preserve the equilibrium variables of the model, which are crucial in the context of flooding simulations.

\section{Summary and outlook} \label{se:summary}

In this paper, we presented a high-order, fully well-balanced, unconditionally positivity-preserving framework for flood simulations.
The discretization based on the notion of production-destruction terms, presented in~\cite{ciallella2022arbitrary}, has been extended to treat general moving equilibria appearing in shallow water systems.
The advantage of this framework lies in the possibility of preserving the positivity of the water height with no constraint on the CFL.
This is a real strength with respect to classical time integration schemes, which experience strong CFL reduction as the order of the method increases, and allows for more realistic applications thanks to the huge computational gain.
In order to achieve the general WB property, while keeping the production-destruction formulation, the high-order reconstruction is blended with a WB one, as proposed in~\cite{berthon2022very,michel2016well}.
This allows to achieve structure preservation for moving equilibria reached after a transient simulation, as shown in \cref{sec:simu_moving}.
On the contrary, when wet-dry simulations are considered, far from existing equilibria, the approach is able to properly perform high-order space and time integration without causing simulation crashes.

There are several perspectives to this work.
They and range from deep questions on the numerical analysis and stability of modified Patankar schemes,
which is an open research topic~\cite{torlo2021stability,izgin2022lyapunov,izgin2022stability}, especially when coupled to space discretizations in the context of PDEs,
to the possible development of this approach on unstructured meshes to exploit advanced mesh adaptation algorithms to capture the flow features with even better resolution, and to save even more computational resources.


\subsection*{Acknowledgements}
M.~C. was funded by a postdoctoral fellowship at ENSAM.
L.~M. was funded by the Schweizerischer Nationalfonds zur F\"orderung der wissenchaftlichen Forschung (SNF) grant 200020\_204917 ``Structure preserving and fast methods for hyperbolic systems of conservation laws'' and by a Postdoc Fellowship at NCSU.
V.~M.-D. acknowledges the support of ANR OptiTrust (ANR-22-CE25-0017).
P.~\"O. was supported by the German Research Foundation (DFG) within SPP 2410, project OE 661/5-1 (525866748) and under the personal grant 520756621 (OE 661/4-1).
D.~T.\ was funded by a SISSA Mathematical Fellowship.
This work has been developed also in the context of the SHARK-FV conference.

All authors would like to thank Jonatan N\'u\~nez for sharing his high-order FV-WENO code on his repository \cite{Nunezrepo}. We have started our work by adapting his code.

\bibliographystyle{abbrv}
\bibliography{sn-bibliography}

\end{document}